\documentclass[a4paper,12pt]{amsart}
\usepackage{amsmath,amssymb,amsthm,amscd,amsbsy,bm}
\usepackage{graphicx}
\usepackage{amsfonts, bbold, bbm}

\usepackage[mathscr]{eucal}

\renewcommand{\a}{\alpha}
\newcommand{\be}{\beta}
\newcommand{\g}{\gamma}
\newcommand{\G}{\Gamma}
\newcommand{\de}{\delta}
\newcommand{\D}{\Delta}
\newcommand{\e}{\epsilon}
\newcommand{\ve}{\varepsilon}
\newcommand{\z}{\zeta}
\newcommand{\y}{\eta}

\newcommand{\ka}{\kappa}
\newcommand{\vk}{\varkappa}
\newcommand{\la}{\lambda}

\newcommand{\m}{\mu}
\newcommand{\n}{\nu}
\newcommand{\x}{\xi}
\newcommand{\X}{\Xi}
\newcommand{\ro}{\rho}

\newcommand{\Si}{\Sigma}

\newcommand{\ta}{\tau}
\newcommand{\f}{\phi}

\newcommand{\vf}{\varphi}

\newcommand{\hi}{\chi}

\newcommand{\om}{\omega}
\newcommand{\Om}{\Omega}

\newcommand{\U}{\Upsilon}

\newcommand{\B}{{\mathbb B}}

\newcommand{\R}{{\mathbb R}}

\newcommand{\Tt}{\mbox{\texttt{T}}}
\newcommand{\St}{\mbox{\texttt{S}}}
\newcommand{\Nt}{\mbox{\texttt{N}}}
\newcommand{\ab}{{\mathbf a}}

\newcommand{\db}{{\mathbf d}}

\newcommand{\gb}{{\mathbf g}}

\newcommand{\nb}{{\mathbf n}}

\newcommand{\rb}{{\mathbf r}}
\newcommand{\sbb}{{\mathbf s}}
\newcommand{\tb}{{\mathbf t}}

\newcommand{\xb}{{\mathbf x}}
\newcommand{\yb}{{\mathbf y}}

\newcommand{\Bb}{{\mathbf B}}

\newcommand{\Cb}{{\mathbf C}}
\newcommand{\Db}{{\mathbf D}}
\newcommand{\Eb}{{\mathbf E}}

\newcommand{\Gb}{{\mathbf G}}
\newcommand{\Hb}{{\mathbf H}}

\newcommand{\Jb}{{\mathbf J}}
\newcommand{\Kb}{{\mathbf K}}
\newcommand{\Lb}{{\mathbf L}}
\newcommand{\Mb}{{\mathbf M}}
\newcommand{\Nb}{{\mathbf N}}

\newcommand{\Qb}{{\mathbf Q}}

\newcommand{\Sbb}{{\mathbf S}}
\newcommand{\Tb}{{\mathbf T}}

\newcommand{\Vb}{{\mathbf V}}

\newcommand{\Zb}{{\mathbf Z}}

\newcommand{\AF}{\mathfrak A}

\newcommand{\BFB}{\mathfrak B}

\newcommand{\dF}{\mathfrak d}
\newcommand{\DF}{\mathfrak D}

\newcommand{\HF}{\mathfrak H}
\newcommand{\hF}{\mathfrak h}

\newcommand{\KF}{\mathfrak K}

\newcommand{\LF}{\mathfrak L}

\newcommand{\pF}{\mathfrak p}

\newcommand{\RF}{\mathfrak R}

\newcommand{\VF}{\mathfrak V}

\newcommand{\UF}{\mathfrak U}

\newcommand{\Ac}{{\mathcal A}}
\newcommand{\Bc}{{\mathcal B}}

\newcommand{\Dc}{{\mathcal D}}
\newcommand{\Ec}{{\mathcal E}}
\newcommand{\Fc}{{\mathcal F}}

\newcommand{\Hc}{{\mathcal H}}

\newcommand{\Lc}{{\mathcal L}}
\renewcommand{\Mc}{{\mathcal M}}

\newcommand{\Nc}{{\mathcal N}}

\newcommand{\Pc}{{\mathcal P}}
\newcommand{\Qc}{{\mathcal Q}}
\newcommand{\Rc}{{\mathcal R}}

\newcommand{\Xc}{{\mathcal X}}
\newcommand{\Yc}{{\mathcal Y}}
\newcommand{\Zc}{{\mathcal Z}}

\newcommand{\dist}{{\rm dist}\,}

\newcommand{\supp}{\hbox{{\rm supp}}\,}

\DeclareMathOperator{\re}{{\rm Re}\,}
\DeclareMathOperator{\rank}{rank}

\newcommand{\diam}{\operatorname{diam\,}}

\newcommand{\codim}{\operatorname{codim\,}}

\newcommand{\loc}{\operatorname{loc}}

\newcommand{\no}[1]{\pmb{|}{#1}\pmb{|}}
\newcommand{\noul}[1]{\pmb{|}{u}\pmb{|}_{l,\hom,{#1}}}
\newcommand{\novl}[1]{\pmb{|}{v}\pmb{|}_{l,\hom,{#1}}}

\newenvironment{dedication}{\itshape\center}{\par\medskip}

\newtheorem{thm}{Theorem}[section]
\newtheorem{cor}[thm]{Corollary}
\newtheorem{lem}[thm]{Lemma}

\newtheorem{proposition}[thm]{Proposition}

\newtheorem{cond}[thm]{Condition}
\theoremstyle{definition}
\newtheorem{defin}[thm]{Definition}

\newtheorem*{rem}{Remark}

\numberwithin{equation}{section}

\begin{document}


\title[Eigenvalues of singular measures]{Eigenvalues of the Birman-Schwinger operator for singular measures: the noncritical case}

\author{Grigori Rozenblum}
\address{Chalmers Univ. of Technology; The Euler International Mathematical Institute, and St.Petersburg State University; NTU Sirius, Sochi, Russia}
\email{grigori@chalmers.se}
\author{Grigory Tashchiyan}
\address{St.Petersburg University for Telecommunications, Russia}

\begin{abstract} In a domain $\Omega\subseteq \mathbb{R}^\mathbf{N}$ we consider compact, Birman-Schwinger type, operators of the form $\mathbf{T}_{P,\mathfrak{A}}=\mathfrak{A}^*P\mathfrak{A}$; here $P$
 is a singular Borel measure in $\Omega$ and $\AF$ is a noncritical  order $-l\ne -\Nb/2$ pseudodifferential operator. For a class of such operators, we obtain estimates and a proper version of H.Weyl's asymptotic law for eigenvalues,
  with order  depending on dimensional characteristics of the measure. A version of the CLR estimate for singular measures is proved. For non-selfadjoint operators of the form $P_2 \mathfrak{A} P_1$ and $\mathfrak{A}_2 P \mathfrak{A}_1$ with singular measures $P,P_1,P_2$ and negative order pseudodifferential operators $\mathfrak{A},\mathfrak{A}_1,\mathfrak{A}_2$ we obtain estimates for singular numbers.
\end{abstract}
\thanks{Acknowledgements. The research of Grigori Rozenblum was supported by the Russian
Science Foundation (Project 20-11-20032).}
\maketitle
\begin{dedication}
In the memory of Mikhail Birman and Mikhail Solomyak, our teachers.
\end{dedication}


\section{Introduction}

In the recent papers  \cite{RIntegral}, \cite{RSh2},  Birman-Schwinger type operators  in a domain $\Om\subseteq\R^\Nb$ were considered, namely, the ones  having the form $\Tb\equiv\Tb_P\equiv\Tb_{P,\AF}=\AF^*P\AF$. Here  $\AF$ is a pseudodifferential operator in $\Om$ of order $-l=-\Nb/2$ and $P=V\m$ is a finite signed Borel measure containing a singular part. For this particular relation between the order of the operator $\AF$ and the dimension of the space,  $-l=-\Nb/2,$ the case we call '\emph{critical}', it was found that the decay rate of the eigenvalues $\la_k^{\pm}(\Tb)$ does not depend on the dimensional characteristics of  measure $\mu$. If measure $\mu$ possesses certain regularity property, namely, it is $s$-regular in the sense of Ahlfors, $s>0$,
the eigenvalue counting function $n_{\pm}(\la,\Tb)=\#\{k: \pm \la^\pm_k(\Tb)>\la\}$ admits estimate $n_{\pm}(\la,\Tb)\le C_{\pm}(V,\m,\AF)\la^{-1}$, with constants depending linearly on the $L\log L$-Orlicz norm of the function $V$ with respect to measure $\m$,  thus the estimate has correct, quasi-classical, order in the coupling constant. This estimate enables us to prove for the eigenvalues of $\Tb$ asymptotic formulas of Weyl type, for the case when $\mu$  is the Hausdorff measure on a certain class of sets, including compact Lipschitz surfaces. Moreover, even in the case when the asymptotics is unknown, the order sharpness of estimates was endorsed by the same order lower estimates.

The reasoning in  \cite{RIntegral}, \cite{RSh2} was based upon the traditional variational approach to the study of eigenvalue distribution, in the version originating in papers by M. Birman and M. Solomyak in late 60-s - early 70-s. In the presence of more modern and more powerful methods, the variational one can still produce new  results involving rather weak regularity requirements in the setting of  spectral problems.

After \cite{RIntegral}, \cite{RSh2}, by our opinion, it is interesting to investigate
  spectral properties for operators of the form $\Tb_P$ with singular measure $P=V\m$ in the \emph{noncritical} case, i.e., for $l\ne\Nb/2$.
Our  results in this paper show that the eigenvalue behavior  in the \emph{subcritical}, $l<\Nb/2$, and in the \emph{supercritical}, $l>\Nb/2$, cases differ essentially from each other, and both from the critical case above.

To explain this difference, we compare the eigenvalue distribution for  an operator with singular measure with the well understood  case of  an absolutely continuous one. In the latter case, $P=VdX$, with rather weak integrability conditions imposed on $V$, for an operator $\AF$ of order $-l$, the eigenvalue counting function has Weyl order (it hints for by the subscript $W$ in the notation), namely, asymptotically
\begin{equation*}
    n_{\pm}(\la,\Tb_{P})\sim C_{\pm}(P,\AF)\la^{-\theta_W},\, \la\to 0, \, \theta_W = \Nb/(2 l).
\end{equation*}
Thus, $\theta_W=1$ in the critical case.
It was established in \cite{RSh2} that for a singular measure, in the critical case, the eigenvalue distribution order does not depend on certain (we call it dimensional) characteristics of  measure $\mu$ and this order still equals $1.$

We recall that a measure $\m$ on $\R^\Nb$ with support (the smallest closed set of full measure) $\Mc$ is called \emph{Ahlfors regular} of order $s$ if for some constants $\Ac=\Ac(\mu),\,\Bc=\Bc(\mu)>0$, the measure of the ball $B(X,r)$ with center $X\in\Mc$ and radius $r$ satisfies
\begin{equation}\label{AR}
    \Bc r^{s}\le \m(B(X,r))\le \Ac r^{s}
 \end{equation}
 for any $X\in \Mc$ and $r\le \diam(\Mc)$. It is known that such measure is equivalent to the dimension $s$  Hausdorff measure $\Hc^s$ on $\Mc$ (see, e.g.,
\cite{DavSem}, Lemma 1.2).
 In the critical case estimates of the eigenvalues of  operator $\Tb_P$ were obtained in \cite{RSh2} under the condition that \eqref{AR} is satisfied for some $s>0;$ moreover, both inequalities in \eqref{AR} were essentially used in the proof, while the decay order of eigenvalues does not depend on $s$.

In the noncritical cases, $2l\lessgtr \Nb$, as we find out in this paper,  the situation is different.   In the \emph{subcritical} case, $2l<\Nb$, previously known results concern mostly  order $-2$ operators ($l=1$) in a domain $\Om$ in $\R^\Nb,$ $\Nb\ge2$ and the singular part of $P$ being supported on a (more or less) regular codimension one surface $\Si\subset\Om;$  here the contribution of such measure to the spectral asymptotics of $n_{\pm}(\la,\Tb_{P})$ is known to have  order $\theta>\theta_W,$ i.e.,  larger than the order of the Weyl term corresponding
 to the absolutely continuous part, see, e.g., \cite{Agr}.

In the present  paper, we establish that
in the subcritical case, such spectral property of a singular measure manifests itself under rather general conditions and the contribution  of such singular measure to the spectrum becomes  stronger as long the dimension of the support of the measure decreases but still satisfies
$s> \Nb-2l$. Moreover, this eigenvalue estimate is established \emph{only} under the sole condition that the \emph{upper} bound in \eqref{AR} holds, namely
\begin{equation}\label{AR+}
    \m(\Mc\cap B(X,r))\le \Ac(\m) r^{s}.
\end{equation}
So, for measures $\mu$ satisfying \eqref{AR+}, we obtain eigenvalue estimates for $n_{\pm}(\la,\Tb)$ with order depending on $\Nb,l$ and the exponent $s\in(\Nb-2l,\Nb)$ in \eqref{AR+},
\begin{equation}\label{n+}
    n_\pm(\la,\Tb)\le C \Ac(\m)^{\theta-1} \int{V_{\pm}}(X)^{\theta}\mu(dX) \la^{-\theta}, \, \theta=\frac{s}{2l-\Nb+s}>1,
\end{equation}
with constant $C$ depending on $\Nb,l,s$ and  operator $\AF.$
Thus,  the order of the eigenvalues decay  and the integrability class of $V$, involved in the eigenvalue  estimate, depend on the exponent $s.$ 


In the usual way, via the Birman-Schwinger principle, estimate \eqref{n+} leads to an estimate of the number $ \nb_-(\DF^*\DF-P)$ of negative eigenvalues of a, properly defined, Schr\"odinger-like operator $\DF^*\DF-P$ with a uniformly elliptic operator $\LF$ of order $l$: under some additional conditions,
\begin{equation*}
    \nb_-(\DF^*\DF-P)\le C (\DF,\m) \int V_+(X)^{\theta}\m(dX).
\end{equation*}
This inequality can be treated as an analogy, for singular measures, of the well-known CLR estimate.

On the other hand, in the \emph{supercritical} case, $2l>\Nb$, the singular part of the measure $\m$ gives to the eigenvalue distribution a \emph{weaker} contribution  than the Weyl term. This kind of spectral behavior of singular measures was, probably, first discovered by
 M. Krein in 1951, see \cite{Krein},
where eigenvalues of the singular string, $-\la u''=Pu,$ on a finite  interval with $P$ being a Borel measure, have been studied. The leading, Weyl order, term in the eigenvalue distribution, according to \cite{Krein}, is determined only by the absolutely continuous component of the measure $P,$ while the singular component makes a weaker contribution.  This spectral problem is equivalent to the one for an operator of the form  $\Tb_{P,\AF}$ with $\Nb=1, l=1$, thus a supercritical one.  Such relative weakness of contribution of the singular part was established later in a  rather general setting  by M.Birman and V.Borzov (\cite{BB}, \cite{Bor}, see also the presentation in \cite{BS}), however, without specifying  the exact order of eigenvalues decay.

In this paper we show that in the supercritical case the \emph{lower} bound in \eqref{AR} only,
namely,
\begin{equation}\label{AR-}
    \m( B(X,r))\ge  \Bc r^{s},\, 0<s\le \Nb,\, \Bc=\Bc(\m),
\end{equation}
is sufficient for the validity for  operator $\Tb_{P,\AF}$ eigenvalue estimates  for a compactly supported  measure $\m,$
\begin{equation}\label{n-}
    n_{\pm}(\la, \Tb_{P,\AF})\le C \la^{-\theta}\Bc^{-\be}\left[\int V_{\pm}(X)\m(dX)\right]^{\theta} \m(\Om)^{1-\theta},
\end{equation}
where (again) $\theta=\frac{s}{2l-\Nb+s},$ but $\theta<1$ now, and $\be =\frac{2l-\Nb}{s}=\theta^{-1}-1.$

 More complicated -- and less sharp -- results are obtained when the compact support condition for $\m$ is dropped. For absolutely continuous measures this problem has been studied earlier, see, e.g.,  \cite{BB}, \cite{BLS}.

An essentially  different approach to obtaining eigenvalue estimates for operators involving singular measures has been elaborated by H.Triebel, \cite{Triebel2}, based upon his earlier studies \cite{Triebel1}, as well as in his  monographs, including the ones joint with D.Edmunds \cite{ET} and joint with D.Haroske \cite{HT}.  These authors consider operators of the form $\Gb=P_2 \BFB(X,D) P_1,$  where $P_j$ are compactly supported functions  on a domain in $\R^\Nb$ (in \cite{ET}, \cite{HT}) or, in the singular case,  measures $V_j\m$ supported on a compact set $\Mc$ whose Hausdorff measure $\m=\Hc^s|_{\Mc}$ satisfies \eqref{AR}, see \cite{Triebel2}, and obtain eigenvalue estimates.  We show that this kind of operators can be dealt with by our method as well, moreover, we obtain more sharp results. We consider also non-selfadjoint operators of the form $\AF_2 P\AF_1$ with pseudodifferential operators $\AF_1,\AF_2$ and $P$ is, as before, a singular measure.
  Additionally, in Appendix, we discuss   approach of \cite{Triebel2} in more detail and  compare the spectral results with ours.

The last part of the paper is devoted to obtaining eigenvalue asymptotics for operators of the form $\Tb_{P,\AF}$ for measure $\m$ being the Hausdorff measure on a Lipschitz surface and, further, on some more general sets of integer dimension, namely, on uniformly rectifiable sets. The general approach to the spectral asymptotics of non-smooth spectral problems, elaborated by M. Birman and M. Solomyak more than 50 years ago, enables us  to derive  asymptotic formulas for the eigenvalues of this kind of compact self-adjoint operators under fairly weak regularity assumptions, as soon as correct order upper eigenvalue estimates involving the integral norm of the weight function are obtained.  Here we follow essentially  the pattern of  \cite{RIntegral} -- \cite{RSh2}, however certain  modifications are needed. In the presentation in this paper, we skip more standard points and concentrate on the essentials. We present also some visual examples demonstrating our general results.

As it has been established in the extensive literature devoted to eigenvalue distribution,
the variational method does not produce more or less sharp values for constants in  eigenvalue estimates. The same shortcoming is present in our paper as well. With few exceptions, we do not care for values of constants in  formulas, however we indicate which important characteristics of a particular  problem these constants depend on. When needed, a constant  is marked by the number of formula where it appears first, say $C_{3.15}$ denotes the constant met first in (3.15). Otherwise, constants are denoted by the symbol $C$ and may change value when passing to the next formula.

Our considerations have their roots in  ideas of the paper \cite{RSh2}, written by the first-named author jointly with Eugene Shargorodsky.   We express deep gratitude to Eugene for benevolent attention and stimulating discussions on the topic.

\section{Estimates in  Sobolev spaces  and singular  Birman-Schwinger operators}

 The starting point in our study is deriving eigenvalue estimates similar to \eqref{n+}, \eqref{n-} for operator $\Sbb_{P,\Om}$ which is defined in the Sobolev space $H^l(\Om)$ or $H^l_0(\Om)$, $\Om\subset \R^\Nb,$ or $\Lb^l$ (see Sect. 2.3 below) by the quadratic form $\sbb_P[v]=\int |v(X)|^2 P(dX)$. Under the conditions we impose on the measure $P$,  form $\sbb_P,$ defined initially on continuous functions in the Sobolev space, admits a bounded extension to the whole  space and thus defines a bounded self-adjoint operator which further on proves to be compact.

\subsection{Sobolev spaces}  In our considerations, $\Om$ is usually  a bounded open set in $\R^\Nb$ with piecewise smooth boundary (of interest are domains with smooth boundary and cubes; such domains are called \emph{nice}). The case $\Om=\R^{\Nb}$ is considered as  well. For  the Sobolev space $H^l(\Om),$  $l>0, $
 the Sobolev norm and the homogeneous Sobolev
 semi-norm are defined by
\begin{align}\label{Sob.norm}
\|u\|^2_{H^{l}(\Om)} := \sum_{|\n|_1 \le l} \left\|\partial^\n u\right\|^2_{L_2(\Om)} , \qquad
\noul{\Om}^2 := \sum_{|\n|_1 = l} \left\|\partial^\n u\right\|^2_{L_2(\Om)},
\end{align}
for an integer  $l $, and by
\begin{align*}
& \|u\|^2_{H^{l}(\Om)} := \|u\|^2_{H^{[l]}(\Om)} +\noul{\Om}^2, \\
& \noul{\Om}^2 := \sum_{|\n|_1 = [l]} \int_\Om\int_\Om \frac{\left|(\partial^\n u)(X) - (\partial^\n u)(Y)\right|^2}{|X - Y|^{\Nb+ 2\tilde{l}}}\, dXdY,
\end{align*}
for a non-integer $l =[l]+\tilde{l}$, $\tilde{l}\in(0,1)$ (here $|\cdot|_1$ denotes the standard norm in $\ell^1$.)

 Consider the case of a bounded domain $\Om$ first. Here, the closure of the space of smooth compactly supported in $\Om$ functions in $H^{l}(\Om)$ is denoted $H^{l}_0(\Om).$ The norm $\no{.}_{l,\hom,\Om}^2$ is equivalent to the standard norm \eqref{Sob.norm} on $H^{l}_0(\Om).$  The important property of the space $H^l(\Om),$ used for obtaining eigenvalue estimates, is the following.
Let $\Om$ be a cube and $H_\top^{l}(\Om)$ denote the subspace in $H^{l}(\Om)$ consisting of functions $L_2(\Om)$-orthogonal to all polynomials of degree less than $l$.
\begin{proposition}\label{EquivNorm}On the subspace $H_\top^{l}(\Om),$ the standard $H^l$-norm is equivalent to its homogeneous part,
\begin{equation}\label{Eq}
  \|u\|^2_{H^{l}(\Om)}\le C(l,\Om)\noul{\Om}^2, \, u\in H_\top^{l}(\Om).
\end{equation}
\end{proposition}
This fact was established in the initial publication by S.L.Sobolev for integer $l$ and became folklore later. An example of a rigorous proof of \eqref{EquivNorm} for a fractional $l$ can be found, e.g., in \cite{SZSol}.

Sobolev spaces for $2l<\Nb$ are closely related to  Riesz potentials,
\begin{equation*}
   I_l(f)(X)=\int_{\Om} |X-Y|^{l-\Nb}f(Y)dY,\, 2l<\Nb;
\end{equation*}
this connection is discussed in details in \cite{MazBook}. For a nice domain $\Om,$ the space of Riesz potentials  $v=I_l(f)$ with $f\in L_2(\Om)$ coincides with the Sobolev space $H^l(\Om)$ and the standard norm in the latter space and the $L_2$ norm of $f$ are equivalent.

\subsection{The case $2l<\Nb$}\label{2l<N.sub}
The basic result we use here  is the generalized Sobolev embedding theorem.
\begin{proposition}\label{Sob.Embed} Let $\Om\subset \R^\Nb$ be a nice bounded domain, $2l<\Nb.$ Suppose that the Borel measure $\mu$ on $\Om$ satisfies \eqref{AR+} with $s>\Nb-2l$ and with constant $\Ac=\Ac(\mu)$.
Then for any $v\in H^l(\Om)\cap C(\Om),$
\begin{equation}\label{Sob}
    \left(\int_{\Om}|v(X)|^q\mu(dX)\right)^{2/q}\le  C_{\ref{Sob}} \Ac^{2/q} \|v\|^2_{H^l(\Om)},\, q=\frac{2s}{\Nb-2l}>2,\, \frac2q=(1-\theta^{-1}),
\end{equation}
with constant $C_{\ref{Sob}}$ not depending on measure $\m$.
\end{proposition}
For integer $l,$ this result is a special case of the generalized Sobolev  theorem, see \cite{MazBook}, Theorem 1.4.5, and for a general  $l$, it is a special case $p=2$ of the D. Adams theorem, see Theorem  1.4.1/2 and Theorem 11.8 in \cite{MazBook}.

For a given complex Borel measure $P=V\mu$ on $\Om$ with $\m$ satisfying \eqref{AR+},  we consider the quadratic form $\sbb[v]=\int |v(X)|^2 P(dX)=\int|v(X)|^2 V(X)\mu(dX)$ defined initially in the space $H^{l}(\Om)\cap C(\Om)$.
We suppose that $V\in L_{\theta,\m},$ $\theta=\frac{s}{2l+s-\Nb}>1.$ By the H\"older inequality and \eqref{Sob},
\begin{gather}\label{H1}
   |\sbb[v]|\le \left[\int_{\Om}|V(X)|^{\theta}\m(dX)\right]^{\frac{1}{\theta}}\left[\int_{\Om}|v(X)|^q\m(dX)\right]^{\frac{2}{q}}\le \\ C_{\ref{Sob}}\nonumber \Ac^{2/q} \left[\int_{\Om}|V(X)|^{\theta}\m(dX)\right]^{\frac{1}{\theta}} \|v\|_{H^l(\Om)}^2.
\end{gather}
It follows that the quadratic form $\sbb[v],$ defined initially on $H^{l}(\Om)\cap C(\Om),$ admits a unique  bounded extension by continuity to the whole of $H^{l}(\Om)$ and thus defines a bounded selfadjoint operator, which we denote by $\Sbb_{l,P,\Om}$ (subscripts or some of them may be dropped if this does not cause confusion.) Further on, we assume this extension already made.

Since on $H^l_0(\Om)$ the Sobolev norm is equivalent to its homogeneous part, it follows from \eqref{H1} that

\begin{equation}\label{H11}
    |\sbb[v]|\le \Zc\novl{\Om}^2, \, v\in H^l_0(\Om),\, \Zc= C_{\ref{H11}} \Ac^{2/q} \left[\int_{\Om}|V(X)|^{\theta}\m(dX)\right]^{\frac{1}{\theta}},
\end{equation}
with constant $C_{\ref{H11}}$ depending only on the domain $\Om, l, \Nb.$

We consider now   operator $\Sbb^{\top}$ defined by the same quadratic form $\sbb[v],$ but on the subspace $H^l_{\top}(\Om)\subset H^l(\Om)$  with norm
  $\novl{\Om}^2,$ in the particular case of $\Om$ being a cube $Q$ in $\R^\Nb.$ We are interested in the dependence of the norm of this operator on the size of the cube $Q.$ If $Q$ is the unit cube $Q_1$, then by \eqref{H1} and Proposition \ref{EquivNorm},
\begin{equation}\label{H1.5}
    \left(\int_{Q_1}|v(X)|^q\m(dX)\right)^{2/q}\le C_{\ref{H1.5}}\Ac^{2/q}\novl{Q_1}^2, \, v\in H^l_\top(Q_1)
\end{equation}
and
\begin{equation}\label{H2}
   |\sbb[v]|\le C_{\ref{H1.5}} \Ac^{2/q}\left[\int_{Q_1}|V(X)|^{\theta}\m(dX)\right]^{\frac{1}{\theta}}\novl{Q_1}^2,
\end{equation}
with certain constant $C_{\ref{H1.5}}$ not depending on  measure $\m.$

Our aim now is to describe the behavior of the constant  in \eqref{H1.5} when the unit cube $Q_1$ is replaced by a cube $Q_t$ with edge $t$. We place the center of co-ordinates at some point in $\Mc$.
 After the scaling  $X\mapsto Y=t^{-1}X,$    cube $Q_t$ is mapped onto $Q_1$.
 Thus, for $v\in H^l(Q_1), \, \tilde{v}(X)=v(t^{-1}X),$ we have $ \tilde{v}\in H^l_\top(Q_t),$ and, by the dilation homogeneity of the norm,
 \begin{equation*}
 \pmb{|}\tilde{v} \pmb{|}_{l,\mathrm{hom},Q_t}^2=t^{2l-\Nb}\novl{Q_1}^2.
\end{equation*}
Recall that the homogeneity order ${2l-\Nb}$ is
 negative in the subcritical case under consideration.
With measure $\m$ on $Q_1$ we associate measure $\tilde{\m}$ on $Q_t$: $\tilde{\m}(E)=t^{s}\m(t^{-1}E),$ $t^{-1}E\subset Q_1.$
Measure $\tilde{\m}(E)$ satisfies
$\tilde{\m}(B(X,r))=t^{s}\m(B(t^{-1}X, t^{-1}r))\le \Ac r^{s}$ with the same constant $\Ac$ as in  \eqref{AR+}.
Therefore, since $s=
\theta(\frac{\Nb}{2}-l),$ after this change of variables, all terms in \eqref{H1.5} acquire the same power of the scaling parameter $t,$ therefore
$|\int_{Q_t}|\tilde{v}(X)|^{q}\tilde{\m(dX)}|^{2/q}\le \Ac^{2/q}C_{\ref{H1.5}} \novl{Q_t}^2, $ $v\in H^{l}_{\top}(Q_t)$
with constant $C_{\ref{H1.5}}$ not depending on the size of the cube $Q_t.$
As a result, by the H\"older inequality, we, finally, obtain
\begin{equation}\label{Hom.sub}
    \left |\int_{Q_t}|v(X)|^2 V(X) \tilde{\mu}(dX)\right|\le C_{\ref{H1.5}} \Ac^{2/q}\left(\int_{Q_t}|V(X)|^{\theta}\tilde{\mu}(dX)\right)^{1/\theta}\novl{Q_t}^2,
\end{equation}
$ v\in H^{l}_{\top}(Q_t).$ We denote $C_0=\max(C_{\ref{H11}},C_{\ref{H1.5}}).$

\subsection{An  estimate in $\R^\Nb$}
For obtaining eigenvalue estimates, we also need a version of the inequality similar to \eqref{H2}, but without the condition of boundedness of $\Om\subset\R^{\Nb},$ including the case of $\Om=\R^\Nb.$ It is sufficient to consider $\Om=\R^\Nb$.
There are equivalent definitions of the space $\Lb^l=\Lb^l(\R^\Nb)$ (called usually the homogeneous Sobolev space or the space of Riesz potentials.) On the one hand, it is the closure of $C_0^{\infty}(\R^\Nb)$ in the metric $\novl{\R^\Nb}^2 =\|(-\Delta)^{l/2} v\|_{L_2}^2$. Equivalently, $\Lb^l$ is the space of Riesz potentials $v(X)=(I_lf)(X)=\int |X-Y|^{2l-\Nb}f(Y)dY,$ $f\in L_2(\R^\Nb)$, with norm equal to the $L_2(\R^{\Nb})$-norm of $f.$ Finally, the space $\Lb^l$ can be described as the space of functions $v\in L_{2,loc}(\R^\Nb)$ such that $(-\Delta)^{l/2} v\in L_2(\R^\Nb)$ in the sense of distributions as well as $|X|^{l-\Nb/2}v(X)\in L_2(\R^\Nb),$ with norm defined by the sum first of these  $L_2$-norms. Locally, the space $\Lb^l$ coincides with the Sobolev space $H^l$, however functions in $\Lb^l$ may have a slower  decay rate at infinity.
\begin{lem}\label{lem.RN}
Let $\m$ be a locally finite Borel measure on $\R^\Nb$ satisfying \eqref{AR+} for $X\in\Mc$, $\Mc=\supp \m,$ and $0<r\le \diam{\Mc}$ (for all $r>0$ if $\Mc$ is unbounded.) Let $V\in L_{\theta,\m},$ $\theta=\frac{s}{2l+s-\Nb}.$ Then
\begin{equation}\label{est.RN}
    \left|\int V(X) |v(X)|^2 \m(dX)\right|\le C_{\ref{H11}} \left(\int |V(X)|^{\theta}\mu(dX)\right)^{\frac1\theta}\Ac^{1-\theta^{-1}}\no{v}_{l,\hom,\R^\Nb}^2, \, v\in \Lb^l.
\end{equation}
\end{lem}
\begin{proof} Since $C_0^\infty(\R^\Nb)$ is dense in $\Lb^l$, it suffices to prove \eqref{est.RN} for $v\in C_0^\infty$. Such embedding of the space of Riesz potentials into $L_{q,\m},$ similar to  inequality \eqref{Sob.Embed},  is, again, contained in \cite{MazBook}, see  there Theorem 11.8;   \eqref{est.RN} follows then by the H\"older inequality.
\end{proof}
Note that unlike the case of the Sobolev space, the inequality for $\m(B(X,r))$ is required to hold for \emph{all} $r>0$ if $\Mc$ is not compact.
\subsection{The case $2l>\Nb$}\label{Sect.cube.2l>N}

We consider now the supercritical case $2l>\Nb.$  Let measure $\m,$ $\supp \mu\subset Q_t,$ satisfy the estimate (converse to \eqref{AR+}):
\begin{equation*}
    \m(B(X,r))\ge \Bc(\m) r^{s}, \, r\le \diam{Q_t}, \, X\in Q_t,
\end{equation*}
for some cube $Q_t$ with edge $t.$
For such cube,  the classical embedding $H^l(Q_t)\subset C(Q_t)$ is valid with estimate

\begin{equation}\label{embedding.C}
    |v(X)|^2\le C t^{2l-\Nb}\novl{Q_t}^2 =C|Q_t|^{2l/\Nb-1} \novl{Q_t}^2,\, X\in Q_t,
\end{equation}
$v\in H^l_\top(Q_t)$, and, further,
\begin{equation}\label{ineq.1}
    \left|\int_{Q_t}V(X)|v(X)|^2\m(dX)\right|\le \int_{Q_t}|V(X)|\m(dX)  \sup_{X\in Q_t}|v(X)|^2.
\end{equation}

It follows from \eqref{AR-} that $t^{2l-\Nb}\le C \Bc^{-\frac{2l-\Nb}{s}}\mu(Q_t)^{\frac{2l-\Nb}{s}}.$
Taking into account \eqref{embedding.C}, \eqref{ineq.1}, we obtain
\begin{equation}\label{Est.cube 2l>N}
\left|\int_{Q_t}V(X)|v(X)|^2\m(dX)\right|\le C_{\ref{Est.cube 2l>N}} \Bc^{-\frac{2l-\Nb}{s}} \m(Q_t)^{\frac{2l-\Nb}{s}} \int_{Q_t}|V(X)|\m(dX) \novl{Q_t}^2,
\end{equation}
$ \, v\in H^l_{\top}(Q_t)$, with constant $C_{\ref{Est.cube 2l>N}}$ not depending on $V, v,t,\mu$.

 Inequality \eqref{embedding.C} is valid for  $v\in H^l_0(Q_t)$    as well  (probably, with a different constant which we still denote by $C_{\ref{Est.cube 2l>N}}$). Therefore estimate \eqref{Est.cube 2l>N} is valid for $v\in H^l_0(Q_t)$.

For a nice bounded domain $\Om,$ a similar estimate holds for \emph{all} functions in $H^l(\Om),$
\begin{equation}\label{Embedding.c.Hl}
    |v(X)|^2\le C\|v\|^2_{H^l(\Om)},
\end{equation}
and
\begin{equation}\label{Norm.est.Hl}
  \left|\int_{\Om}V(X)|v(X)|^2\m(dX)\right|\le C \int_{\Om}|V(x)|\m(dX) \|v\|^2_{H^l(\Om)},
\end{equation}
but with \emph{certain} dependence of the constant $C$ in \eqref{Embedding.c.Hl}, \eqref{Norm.est.Hl} on the domain $\Om$. We, however, will use it only for $\Om$ being a unit cube, and here the constant is absolute.

\section{Eigenvalue estimates for $\Sbb_{V\m},$ and the CLR bound for singular measures}\label{Sect.3}
\subsection{Preparation}\label{prep} We will use the following geometrical statement established in \cite{RSh2} (the two-dimensional version was proved in \cite{KarSh}).
 For a fixed $\Nb$-dimensional cube $\Qb\subset \R^\Nb$ we call  cube $Q$ parallel to $\Qb$ iff all (one-dimensional) edges of $Q$ are parallel to the ones of $\Qb.$
 \begin{lem}\label{Geom.Lemma} Let $\m$ be a locally finite Borel measure in $\R^\Nb$ not containing point masses. Then there exists a cube $\Qb\subset\R^{\Nb}$ such that for any cube $Q$ parallel to $\Qb$, $\m(\partial Q)=0.$
\end{lem}
Another important ingredient is the Besicovitch covering lemma, see, e.g., Theorem 1.1 in \cite{Gu}.
\begin{lem}\label{Besik} Let $\Mc$ be a set in  $\R^{\Nb}$ and with each point $X\in \Mc$ an open cube $Q_X$ centered at $X$ is associated, all cubes are parallel to each other,  and the size of cubes is uniformly upper bounded. Then there exists a covering $\U$ of $\Mc$ consisting of cubes $Q_X$ so that $\U$ can be  split into a controllably finite collection of families $\U_j,\, j=1,\dots,\ka=\ka(\Nb) $ such that cubes in each of $\U_j$ are disjoint.
\end{lem}
\subsection{Estimates in the subcritical case $2l<\Nb$}\label{2l.l.N}

\begin{thm}\label{Thm.Est.2l<N}
Let $\Om\subset\R^\Nb$ be a bounded open set and $\mu$ be a measure on $\Om$ satisfying \eqref{AR+} with some $s>\Nb-2l.$ Let $V$ be a real function on $\Mc=\supp\m$, $V\in L_\theta\equiv L_{\theta,\m},$  $\theta=\frac{s}{s+2l-\Nb}.$
Then for  operator $\Sbb=\Sbb_{l,V\m,\Om}$ defined by the quadratic form $\sbb[v]=\int V|v(X)|^2\mu(dX)$ in the Sobolev space $H^l_0(\Om),$ the following  estimate holds:
\begin{equation}\label{BS+}
    n_\pm(\la,\Sbb)\le C_{\ref{BS+}}\Ac^{\theta-1}\int V_{\pm}(X)^\theta\mu(dX)\la^{-\theta},
\end{equation}
with constant $C_{\ref{BS+}} =C_{\ref{BS+}}(l,s,\Nb).$\\
If $V$ is a complex-valued function,  $V\in L_\theta$, then the estimate, similar to \eqref{BS+}, holds for the distribution function  $n(\la,\Sbb)$ of the singular numbers of $\Sbb:$
\begin{equation}\label{BSsing}
    n(\la,\Sbb)\le C_{\ref{BSsing}}\Ac^{\theta-1}\int |V(X)|^\theta\mu(dX)\la^{-\theta}
\end{equation}
\end{thm}\begin{proof}  Since any complex function  is a linear combination of $4$ nonnegative ones,  it is sufficient to consider the case of a nonnegative density $V$ and prove the estimate for $n_{+}(\la,\Sbb).$
Denote $\Kb=\Zc^{\theta}\equiv C_{\ref{H11}}^{\theta}\Ac^{\theta-1}\int V(X)^\theta\mu(dX).$
The structure of the proof is the following. We find constants, $\n_1$ and $\n_2$ depending only on $l,\Nb,s$ such that for $\la > \n_1\Kb^{\frac1\theta},$ we have $n_+(\la,\Sbb)=0,$ while for $\la < \n_2\Kb^{\frac1\theta}$ the required estimate holds, $n_+(\la,\Sbb)\le C \Kb \la^{-\theta}$. From these inequalities, it follows that \eqref{BS+} holds for \emph{all} $\la>0,$ probably, with a different constant $C$. In fact, we need to establish our estimate only for $\la\in [ \n_2\Kb^{\frac1\theta}, \n_1\Kb^{\frac1\theta}].$ For $\la$ in this interval, the number $\la \n_2 \n_1^{-1}$ is smaller than $\la$ and smaller than $\n_2\Kb^{\frac1\theta},$ and therefore,  for such $\la:$

\begin{equation*}
    n_+(\la,\Sbb)\le n_+(\la \n_2 \n_1^{-1},\Sbb)\le C  (\la \n_2 \n_1^{-1})^{-\theta}\Kb= C (\n_2 \n_1^{-1})^{-\theta}\la^{-\theta}\Kb,
\end{equation*}
exactly what we needed.

Due to the obvious homogeneity,
\begin{equation*}
    n_{\pm}(\la,\Sbb)=n_{\pm}(\g\la,\g\Sbb), \,\g>0,
\end{equation*}
it is sufficient to consider the case  $\Kb=1.$

Following the way of reasoning described above, we consider large $\la$ first. By \eqref{H11}, the norm of operator $\Sbb$ is not greater than $\Kb=1$,
 therefore,
 $n_+(1,\Sbb)=0,$ so we can take $\n_1=1.$

Next, for sufficiently small $\la<\n_2$,
following the variational principle, we will construct a subspace $\Lc(\la)$ of codimension not greater than $C\la^{-\theta}$ such that
$\sbb[v]<  \la  \novl{\Om}^2,$ $v\in\Lc(\la), \, v\ne 0.$

We act by an adaptation of the original construction by M. Birman-M. Solomyak, with contribution by G.Rozenblum, see e.g., \cite{BirSolLect}. Consider the function of $\m$-measurable sets,

\begin{equation*}
\Jb(E)=2C_{0}\Ac^{\theta-1}\int_{E}V(X)^{\theta}\mu(dX)
\end{equation*}
(measure $\mu$ is supposed to be extended by zero to $\Om\setminus\Mc.$) Thus, $\Jb(\Mc)=2.$ We fix a cube $\Qb\subset\R^\Nb$ according to Proposition \ref{Geom.Lemma}.

Having $\la$ fixed, for each point $Y\in\Mc,$ we consider the family $Q_\ta(Y)$ of cubes parallel to $\Qb$, centered at $Y$ and with edge $\ta$. It follows from
  \eqref{Geom.Lemma} that $\Jb(Q_\ta(Y))$ is a monotonous continuous function of $\ta$ variable; it equals zero for $\ta=0$ and equals $2$ for large
   $\ta$. Due to this continuity, for a given integer $N>0$, to be determined later, there exists a value $\ta=\ta(Y)$ (not necessarily unique) such that $\Jb(Q_{\ta(Y)}(Y))=N^{-1}$.
    Such cubes $Q_{\ta(Y)}(Y)$  form a covering of $\Mc$. Let $\U$ be the subcovering of $\Mc$ consisting of cubes $Q_{\ta(Y)}(Y)$, found according to the Besicovitch
    lemma, $\U=\cup_{1\le j\le \ka} \U_j$. It follows, in particular, that the multiplicity of the covering $\U$ is not greater than $\ka.$
By the disjointness property and additivity of $\Jb$, for each fixed $j,$ we have
\begin{equation}\label{union}
|\U_j|N^{-1}=\sum_{Q\in\U_j}\Jb(Q)= \Jb(\cup_{Q\in\U_j} Q)\le \Jb(\Mc)=2,
\end{equation}
 Therefore, the quantity $|\U_j|$ of cubes in $\U_j$ is not greater than $2N,$ and altogether,
  the quantity of cubes in $\U$ is not greater than $2\ka N.$

Using this covering, we construct the subspace $\Lc(\la).$ Consider the linear space $\Pc(l,\Nb)$ of (real) polynomials in $\R^\Nb$
 having degree less than $l$;
we denote its dimension by $\pF\equiv\pF(l,\Nb)$. Now we consider the collection of linear functionals $\Psi(l,\Nb,\U)$
 on $H^l_0(\Om)$ of the form $\psi(v)=\psi_{p,Q}(v)=\int_{Q}v(X){p(X)}dX$, where $Q$ are cubes in $\U$ and $p$ is in
  a fixed basis in $\Pc(l,\Nb)$; these functionals are, obviously, continuous on $H^l(\Om)$. Thus, there are altogether $\pF|\U|\le 2\ka\pF N$
   functionals in $\Psi(l,\Nb,\U).$ We denote by $\Lc[N]$ the intersection of null spaces of these functionals in $H^l_0(\Om)$. It is a subspace
    of codimension not greater than $ 2\pF\ka N$. This will be the space $\Lc(\la)$ we are looking for, with proper relation of $\la$ and $N$ to be determined.

To find this relation we estimate the quantity $\int |v(X)|^2 V(X)\mu(dX)$ for $v\in\Lc[N]$. Since cubes $Q\in\U$ form a covering of $\Mc=\supp(\m),$
\begin{equation}\label{est1}
  \sbb[v]\equiv  \int_\Mc|v(X)|^2V(X)\m(dX)\le\sum_{Q\in\U}\int_Q|v(X)|^2V(X)\m(dX).
\end{equation}
We recall that for each  cube $Q\in\U$, the restriction of function $v$ to $Q$ belongs to $H^l_{\top}(Q),$ and, therefore, \eqref{est1} and \eqref{Hom.sub} imply
\begin{equation}\label{est3.30}
     \sbb[v]\le
     C\sum_{Q\in\U}\novl{Q}^2\Jb(Q)^{\frac{1}{\theta}}
\le C  \sum_{Q\in\U}\novl{Q}^2 N^{-\frac{1}{\theta}}.
\end{equation}
From the definition of the homogeneous norm $\pmb{|.|}^2_{l,\hom}$ and the finite multiplicity of the covering
 $\U,$ it follows that  $\sum_{Q\in\U}\novl{Q}^2\le \ka\novl{\Om}^2$. Therefore, by \eqref{est3.30},
\begin{equation}\label{est3.1}
    \sbb[v]\le  C_{\ref{est3.1}} \ka  N^{-\frac{1}{\theta}}\novl{\Om}^2.
\end{equation}
Now, for a given $\la>0,$ we take an integer $N=N(\la)$ so that $\la>C_{\ref{est3.1}} \ka  N(\la)^{-\frac{1}{\theta}},$ i.e.,
\begin{equation}\label{N(l)}
N(\la)>(C_{\ref{est3.1}} \ka)^{\theta}\la^{-\theta}.
\end{equation}
 Then \eqref{est3.30} implies
\begin{equation*}
    \sbb[v]\le \la\novl{\Om}^2,\, v\in\Lc(\la),
\end{equation*}
as we wish. In order to estimate the codimension of the subspace $\Lc(\la)$, we need a converse inequality between $\la$ and $N$. Suppose that
\begin{equation}\label{NL}
  (C_{\ref{est3.1}}  \ka)\la^{-1}>1.
\end{equation}
Then  such integer $N=N(\la),$ satisfying   \eqref{N(l)} and $N(\la)<2C_{\ref{est3.1}}^{\theta} \ka^{\theta}\la^{-\theta},$  exists. Due to the last
inequality, the codimension of the subspace $\Lc(\la)$ satisfies
\begin{equation*}
    \codim(\Lc(\la))\le 2\pF\ka N(\la)\le 4 \pF \ka C_{\ref{est3.1}} \ka^{\theta}\la^{-\theta}.
\end{equation*}
Therefore, to find the required $N,$ we need  condition \eqref{NL} be fulfilled, i.e., $\la<\n_2$, where $\n_2=(4\pF C_{\ref{est3.30}})^{\frac1\theta}\ka .$
So, this number $\n_2$ fits into the structure explained in the beginning of the proof.
\end{proof}

A similar estimate, but without control of the dependence of constant on the domain $\Om$ is valid for the operator $\Sbb^{\Nc}$ defined by the same quadratic form $\sbb[v]$, but on the Sobolev space $H^l(\Om),$ a generalization of the operator with Neumann conditions.
\begin{thm}\label{Th.Est.N} Let $\Om\subset\R^\Nb$ be a nice bounded domain, $V\in L_{\theta,\m}$. Then the eigenvalues of operator $\Sbb^{\Nc}_{l,P,\Om}$ satisfy
\begin{equation*}
n_{\pm}(\la,\Sbb^{\Nc})\le C \Ac^{\theta-1}\int V_{\pm}(X)^{\theta}\m(dX)\la^{-\theta}.
\end{equation*}
\end{thm}
\begin{proof}We use the bounded  extension operator $\Ec: H^l(\Om)\to H^{l}_0(\tilde{\Om}),$ $\|\Ec\|=C(\Om)$ for some bounded domain $\tilde{\Om}\supset \Om.$ By the variation principle, $n_\pm(\la, \Sbb^\Nc_{l,P,\Om})\le n_\pm(\la C(\Om)^{-1}, \Sbb_{l,P,\tilde{\Om}})$ and \eqref{BS+} applies.
\end{proof}

\subsection{The case $2l<\Nb.$ The estimate in $\R^\Nb$ and the CLR-type estimate}\label{RLC}
A standard trick enables us to carry over the eigenvalue estimates in Theorem \ref{Thm.Est.2l<N} to the case of operator acting in the whole space. The first stage consists in extending the results to operators  in $\R^\Nb$ but with a measure having compact support.
\begin{proposition}\label{Prop.comp.sup}Let $\Om=\R^\Nb$, $2l<\Nb$; let $\m$ be a measure on $\R^\Nb$ with compact support, satisfying \eqref{AR+}  with $s> N-2l$, and $V\in L_{\theta,\m},$ $\theta =\frac{s}{2l-\Nb+s}>1.$ Then for the operator $\Sbb=\Sbb_{P,\R^\Nb}$  in $\Lb^l,$ the estimate holds
\begin{equation}\label{EstRN.1}
    n_{\pm}(\la,\Sbb)\le C \Ac^{\theta-1}\int_{\R^\Nb}V_{\pm}(X)^{\theta}\m(dX)\la^{-\theta},\, .
\end{equation}
with constant not depending on the size of support of  measure $\m.$
\end{proposition}
\begin{proof} As usual, by the min-max principle, it suffices to consider the case of a positive density $V$ and study $n_{+}(\la,\Sbb).$ Let, for some $\la>0,$
$n_+(\la,\Sbb_{P,\R^{\Nb}})> \nb$ for a certain  integer $\nb$. Then, due to the variational principle and the density of $C_0^\infty(\R^\Nb)$ in $\Lb^l,$ there exists a subspace $\Lc_\nb\subset C_0^\infty(\R^\Nb)$ of dimension $\nb$ such that $\sbb[v]>\la \novl{\R^\Nb}^2, v\in \Lc_\nb\setminus{\{0\}}.$ However, since the space $\Lc_\nb$ is finite-dimensional, there exists a common compact support, i.e., all functions in $\Lc_\nb$ have support in a certain ball $\Om'$ containing the support of measure $\m$. Thus we obtain a subspace of dimension $\nb$ in $C_0^\infty(\Om'),$ where $\sbb[v]>\la \novl{\Om}^2$. However, by Theorem \ref{Thm.Est.2l<N}, this subspace may not have dimension greater than $C \la^{-\theta}\Ac^{\frac{2\theta}{q}}\int V(X)^{\theta}\m(dX)$. This last observation gives us the required inequality $\nb\le C \la^{-\theta}\Ac^{\frac{2\theta}{q}}\int V(X)^{\theta}\m(dX).$
\end{proof}
Now we dispose of the condition of compactness of the support of measure $\mu.$
\begin{thm}\label{ThmRN} Let $\Om=\R^\Nb$, $2l<\Nb$, $\m$ be a measure on $\R^\Nb$  satisfying \eqref{AR+}  with $s> N-2l$, and $V\in L_{\theta,\m},$ $\theta =\frac{s}{2l-\Nb+s}>1.$
Then  estimate \eqref{EstRN.1} holds. Similarly to \eqref{BSsing}, an estimate of the same form holds for singular numbers of operator $\Sbb$ for a complex-valued function $V.$
\end{thm}

\begin{proof}As usual, we may suppose that $V\ge0.$ For a fixed $\la>0,$ we split $V$ into two parts, $V_\la$ and $V_\la',$ so that $V_\la$ has compact support and $V_\la'$ is small, $ C_0 \Ac^{\frac{2\theta}{q}}\int(V_\la')^{\theta}\mu(dX)<(\la/2)^{\theta}.$ Correspondingly, operator $\Sbb$ splits into the sum of two operators, $\Sbb=\Sbb_\la+\Sbb_{\la}'.$ For $\Sbb_\la$ we have eigenvalue estimate of the form \eqref{EstRN.1}
with $V_{\la}$ replacing $V$, and for $\Sbb_{\la}'$ we have the norm estimate, by Lemma \ref{lem.RN}, $\|\Sbb_{\la}'\|< \la/2$, which is equivalent to:
\begin{equation*}
    n_+(\la/2, \Sbb_{\la}')=0,
\end{equation*}
Now, by the Ky Fan inequality,
\begin{equation}\label{KyFan}
    n_+(\la,\Sbb)\le n_+(\la/2, \Sbb_{\la})+n_+(\la/2, \Sbb_{\la}').
\end{equation}
The second term in \eqref{KyFan} equals zero, and Proposition \ref{Prop.comp.sup} applied to $ \Sbb_{\la}'$ gives us the required estimate.
\end{proof}
\begin{rem}\label{RemRN}Since $(-\D)^{-l/2}$ is an isometric isomorphism of $L_2(\R^\Nb)$ onto $\Lb^l,$
the last result can be expressed as an eigenvalue estimate for operator
\begin{equation*}
\Tb_{P,(-\D)^{-l/2}}=(-\D)^{-l/2}P(-\D)^{-l/2}
\end{equation*}
in $L_2(\R^\Nb)$:
\begin{equation*}
    n_{\pm}(\la,\Tb_{P,(-\D)^{-l/2}})\le C \Ac(\m)^{\theta-1}\int_{\R^\Nb}V_{\pm}(X)^{\theta}\m(dX)\la^{-\theta}.
\end{equation*}
\end{rem}

As an automatic consequence of Theorem \ref{ThmRN} we obtain a version of the CLR estimate for singular measures.
\begin{thm}\label{Thm.CLR}Let $2l<\Nb$, let measure $\m$ satisfy \eqref{AR+} with some $s>\Nb-2l,$ and density $V\ge0$ satisfy $\int V(X)^{\theta}\mu(dX)<\infty,$ $\theta=\frac{s}{2l-\Nb+s}>1$. Consider the Schr\"odinger-type operator $\HF=\HF(l,P)=(-\Delta)^l-P$ defined in $L_2(\R^\Nb)$ by means of the quadratic form $\hF[v]=\|(-\Delta)^{l/2}v\|^2_{L_2}-\int|v(X)|^2 P(dX),$ $P=V\m.$ Then for  $N_-(\HF),$ the number of negative eigenvalues of $\HF,$ the estimate holds

\begin{equation}\label{CLR}
    N_-(\HF)\le C(\Nb,l)\Ac(\m)^{\theta-1}\int V(X)^\theta\mu(dX).
\end{equation}
\end{thm}
\begin{proof} The proof, actually, the derivation of \eqref{CLR} from  estimate \eqref{EstRN.1}, is a quite standard application of the Birman-Schwinger principle. We repeat it for the sake of completeness. The quantity $N_-(\HF)$ equals the minimal codimension of subspaces $\Lc\subset C_0^\infty(\R_\Nb),$ for which
 \begin{equation*}
 \int |v(X)|^2V(X)\mu(dX)\ge \novl{\R^{\Nb}}^2,\, v\in \Lc\setminus\{0\}.
 \end{equation*}

 But this codimension is exactly $n_+(1,\Sbb_{P,\R^\Nb}),$ and for the latter quantity we already have the required estimate.
\end{proof}
Later, we obtain a more general version of \eqref{CLR}.

\subsection{Eigenvalue estimates for $\Sbb, \,2l>\Nb$}\label{sect.2lGN}
The aim of this subsection is to prove  eigenvalue estimates in the supercritical case $2l>\Nb.$ The estimates have a somewhat different form, compared with the subcritical case, but the proof is quite analogous.
\begin{thm}\label{Est.2l>N.bdd}Let $\Om
\subset \R^\Nb$ be a nice bounded domain, $\mu $ be a compactly supported finite Borel measure on $\Om$ satisfying \eqref{AR-} with $s>0$ and $V$ be a real-valued function in $L_1=L_{1,\mu}.$
Then for the operator $\Sbb=\Sbb_{l,V\m,\Om}$ defined by the quadratic form $\sbb[v]$ in $H^l_0(\Om)$ the eigenvalue estimate holds
\begin{gather}\label{BS-}
n_\pm(\la,\Sbb)\le C  \Kb\la^{-\theta}, \Kb= \Bc^{-\be} \left(\int V_{\pm}(X)\mu(dX)\right)^{\theta}\m(\Om)^{1-\theta},\\\nonumber
 \be=\theta^{-1}-1=\frac{2l-\Nb}{s}>0, \theta=\frac{s}{2l-\Nb+s}<1,
\end{gather}
with constant $C$ not depending on  measure $\mu$ and  weight function $V.$
The same estimate holds for a nice bounded domain $\Om$ for operator $\Sbb^{\Nc}$  defined by quadratic form $\sbb[v]$ in the Sobolev space $H^l(\Om), $ now with constant depending on $\Om.$
\end{thm}

Note that, unlike the subcritical case, measures possessing point masses are not excluded.
\begin{proof} The statement about the operator in the space $H^l(\Om) $ is reduced to the one for the space $H^l_0(\Om),$ using the extension operator, as it was done for $2l<\Nb$ in Theorem  \ref{Thm.Est.2l<N}. So, we consider operator $\Sbb_{l,P,\Om}$ in the space  $H^l_0(\Om);$ the reasoning here is analogous to the one in the latter Theorem, with the same separate treatment of larger and smaller $\la.$
We explain the modifications in detail, when needed. Again, it suffices to consider only the case $V\ge0$, and, by homogeneity, we can set $\Kb=1.$

 The estimate for large $\la$
follows from the inequality \eqref{Est.cube 2l>N}. For small $\la$, the reasoning, as in the case $2l<\Nb$, consists in constructing a subspace $\Lc(\la)$ of controlled codimension so that $\sbb[v]\le \la \|v\|^2_{H^l_0(\Om)},$ $v\in\Lc(\la).$ We introduce a function of $\m$-measurable sets:
 $$\Jb(E)=C_{\ref{Est.cube 2l>N}} \Bc^{-\be}\left[\mu(E)^{\be}\int_{E}V(X)\m(dX)\right]^{\frac{1}{\be+1}}=
 C_{\ref{Est.cube 2l>N}} \Bc^{-\be}\left(\int_{E}V(x)dx\right)^{\theta}\mu(E)^{1-\theta}.$$ By Proposition \ref{Geom.Lemma}, for each point $Y\in\Mc,$ for a family of concentric cubes $Q_\ta(Y)$ parallel to a certain fixed cube $\Qb$  and centered at $Y\in\Mc$ with edge $\ta$,  $\Jb(Q_\ta(Y))$ is a continuous nondecreasing  function of $\ta$ variable. Therefore, for an integer $N,$ to be determined later, there exists $\ta(Y)$ such that $\Jb(Q_{\ta(Y)}(Y))=(2N)^{-1}.$  From the covering $\{Q_{\ta(Y)}(Y),\, Y\in \Mc\},$ by the same Besicovitch covering lemma, we can extract a finite subcovering $\U$ which can be split into no more than $\ka=\ka(\Nb)$ families $\U_j,$ each consisting of disjoint cubes. Next, we evaluate the number of cubes in $\U.$
For a fixed $j,$ we have, by the H\"older inequality,

\begin{gather*}
(2N)^{-1} |\U_j|=
\sum_{Q\in \U_j}\Jb(Q)= C_{\ref{Est.cube 2l>N}}\Bc^{-\be}\sum_{Q\in\U_j}\left[\int_{Q}V(X)\m(dX)\right]^{\theta}\mu(Q)^{1-\theta}\le\\ \nonumber
 C_{\ref{Est.cube 2l>N}}\Bc^{-\be}\left[\sum_{Q\in \U_j} \int_{Q}V(X)\m(dX)\right]^{\theta}\left[\sum_{Q\in \U_j}\m(Q) \right]^{1-\theta}\\\nonumber
\le C_{\ref{Est.cube 2l>N}}\Bc^{-\be} \left[\int_{\Om}V(X)\m(dX)\right]^{\theta}\m(\Om)^{1-\theta}=\Jb(\Om)=1.
\end{gather*}
It follows that  $|\U_j|$   is not greater than  $2N$ and the quantity of cubes in the whole covering $\U$ satisfies
\begin{equation*}
   |\U|\le 2\ka N.
\end{equation*}

With each cube $Q\in\U$ we associate a collection of linearly independent functionals, scalar products in $L_2(Q)$ with polynomials of degree
 less than $l$; as before, there are $\pF=\pF(\Nb,l)$ of them.
As in Theorem  \ref{Thm.Est.2l<N}, we define $\Lc[N]$ as the common null-space of these functionals. Its codimension is
not greater than $\pF |\U|\le 2\pF\ka N.$ Next, similar to \eqref{est3.1},
we estimate $\sbb[v]$. Due to the finite multiplicity of the covering $\U$ and \eqref{Est.cube 2l>N}, we have
\begin{gather}\label{collection2l>N}
    \sbb[v]\le \ka \sum_{Q\in \U}\int_{Q}|v(X)|^2V(X)\m(dX)\le \\\nonumber
    C \Bc(\m)^{-\be}\sum_{Q\in \U}\m(Q)^{\be}\int_{Q}V(X)\m(dX)\novl{Q}^2=C \sum_Q\Jb(Q)^{\be+1}\novl{Q}^2\\\nonumber \le
    C_{\ref{collection2l>N}}N^{-(\be+1)}\novl{\Om}^2.
\end{gather}
Finally, we take $N> C_{\ref{collection2l>N}}^{\theta}\la^{-\theta}$, so that, by \eqref{collection2l>N}, $\sbb[v]\le \la \novl{\Om}^2;$
 $\theta =(\be+1)^{-1},$  and we can see that for $\la$ small enough, such integer $N$ can be chosen, simultaneously, smaller than
 $2C_{\ref{collection2l>N}}^{\theta}\la^{-\theta}$, which gives for the subspace $\Lc(\la)=\Lc[N]$ the required estimate for its codimension:
 $\codim(\Lc(\la))\le C \la^{-\theta}$.
 \end{proof}

\subsection{The case $2l>\Nb$: the  Birman-Borzov estimate in the whole space}\label{Sect.BBest}
For a measure $\m$ with noncompact support, it follows from \eqref{AR-} that $\mu(\Mc)=\infty$,
 and inequality \eqref{BS+} becomes trivial, i.e.,  does not give any estimate for eigenvalues of operator $\Sbb$. Since the paper \cite{BB} by M.Birman and V.Borzov, quite a lot of work has been done to obtain eigenvalue estimates for this operator as well as for the Schr\"odinger-like operator $\HF=(-\Delta)^l-P$ in $\R^\Nb,$ $2l\ge\Nb$ for an absolutely continuous measure $P,$ see, especially, \cite{BLS} and \cite{Shar2D}. We will not discuss all possible versions of these results for singular measures, and restrict ourselves to just several typical ones.
 An estimate  in \cite{BB}, probably, the first one obtained for the case $2l>\Nb$ in the whole space,  bounds the eigenvalues of the problem $\la[(-\Delta)^l+1] u=Vu$ by a sum of powers of  $L_1$- norms of the function $V$ over the lattice of unit cubes. Later this result was considerably generalized in many directions, however a common feature remained: for operators with fast eigenvalue decay, such as $n_{\pm}(\la)=O(\la^{-\theta}),\, \theta<1,$ estimates should involve the sum of $L_1$- norms, taken to some power, of the weight function over a proper system of compact sets.  To carry over all these results to the case of singular measures considered in this paper  would be a huge task. We demonstrate such generalization of the initial Birman-Borzov result.

\begin{thm}\label{Ex.BB}Let $\mu$ be a measure on $\R^\Nb,$ with noncompact support $\Mc=\supp(\mu)$, satisfying a local version of  condition \eqref{AR-}, namely
\begin{equation}\label{AR-noncomp}
    \m(B(X,r))\ge\Bc(X)r^s, \, X\in \Mc, 0<r<\sqrt{\Nb}, \Bc(X)>0
\end{equation}
(recall, $\sqrt{\Nb}$ is the diameter of the unit cube in $\R^{\Nb}$.) Consider the lattice $\LF$ of unit cubes in $\R^\Nb,$ denote by $\LF_\m$ the set of those cubes $Q$ in $\LF$ that satisfy $\m(Q)>0$  and by  $\Bc(Q), Q\in \LF_\m $ we denote the quantity $\Bc(Q)=\inf_{X\in Q\cap\Mc}\Bc(X).$
Suppose that the density $V\in L_{1,\loc,\mu}$ satisfies 
\begin{equation}\label{BBcond}
    \Mb(\m,V):=\sum_{Q\in\LF_\mu}\Bc(Q)^{1-\theta^{-1}}\left[\int_{Q}|V(X)| \m(dX)^{\theta}\right]\m(Q)^{1-\theta}<\infty
\end{equation}
for $\theta=\frac{s}{2l-\Nb+s}.$ Consider operator $\Sbb=\Sbb_{l,V\m}$ defined by the quadratic form $\sbb[v]$ in the Sobolev space $H^l(\R^\Nb)$, $2l>\Nb.$
Then this operator is bounded, compact, and satisfies
\begin{equation}\label{BBest}
    n_{\pm}(\la,\Sbb)\le C \Mb(\m,V_{\pm})\la^{-\theta}.
\end{equation}
\end{thm}
\begin{proof}Having our estimate \eqref{BS-}, the proof of  \eqref{BBest} can be constructed by repeating the reasoning in \cite{BB}, where an eigenvalue estimate in the whole space was obtained by summing estimates in separate cubes. In the present understanding, this proof fits in just a few lines.
Consider, for each $Q\in \LF_\mu$, operator $\Sbb_Q$ which is defined by the quadratic form $\sbb_Q[v]=\int_{Q}|v(X)|^2 V(X)\m(dX)$ in $H^l(Q)$. The eigenvalues can only grow if we replace the base space $H^l(\R^\Nb)$ by $\oplus H^l(Q)$ with the sum over $Q$ in the lattice $\LF.$ Operators $\Sbb_Q$ act in orthogonal subspaces $H^l(Q)$, therefore,  by the variational principle,
$n_{\pm}(\la, \Sbb)\le \sum_{Q\in\LF_\mu}n_{\pm}(\la,\Sbb_Q),$ and  by summing estimate \eqref{BS-} over cubes $Q$, we arrive at \eqref{BBest}.\end{proof}
Now, since operator $(1-\Delta)^{-l/2}$ is an isometric isomorphism of $L_2(\R^\Nb)$ onto $H^{l}(\R^\Nb),$ Theorem \ref{Ex.BB}  leads to an eigenvalue estimate for a singular Birman-Schwinger operator.
\begin{cor}\label{BB.cor}Let measure $\m$ and density $V$ satisfy conditions of Theorem \ref{Ex.BB}.
Then for the operator $\Tb=(1-\Delta)^{-l/2}(V\m)(1-\Delta)^{-l/2}$ the eigenvalue estimate \eqref{BBest} holds.
\end{cor}

\section{Pseudodifferential Birman-Schwinger operators with singular weight. Spectral estimates}\label{Sect.4}
In this section we  give a detailed  definition of pseudodifferential operators with singular weight and study spectral estimates for the noncritical case(s). We mostly follow the presentation in \cite{RSh2}, \cite{RIntegral}, where the critical case was considered, and emphasize only essential differences.

\subsection{Operators in bounded domains} Let $\Om$ be a bounded open set in $\R^\Nb$ and $\AF$ be a classical pseudodifferential operator of order $-l\ne -\Nb/2,$ with principal symbol $\ab_{-l}(X,\X).$ We suppose that $\AF$ is an operator with compact support in $\Om$ which means that it contains cut-offs to some proper subdomain $\Om'\Subset\Om$,
$\AF=\chi \AF \chi,$ $\chi\in C^\infty_0(\Om).$

Such operator $\AF$ maps $L_2(\Om)$ to $H^l_0(\Om)$, its essential norm in these spaces is bounded by $C(\Nb,l) \sup_{(X,\X) \in \St^*(\Om)}|\ab_{-l}(X,\X)|,$ where $\St^*(\Om)$ is the cospheric bundle of $\Om$.

Let $\m$ be  a singular measure  supported in $\Om$ and $V$ be a real $\mu$-measurable function,  $P=V\mu$.
For $u\in L_2(\Om)$ we consider  the quadratic form $\tb[u]=\tb_{P,\AF}[u]=\tb_{P,\AF,\Om}[u]$ defined as
\begin{equation}\label{Pseu.1}
    \tb_{P,\AF,\Om}[u]=\int_{\Om}|(\AF u)(X)|^2 P(dX).
\end{equation}

Under the conditions we impose now on  measure $\m$ and density $V,$ it will be shown that the quadratic form \eqref{Pseu.1} is well-defined on $L_2(\Om)$ and determines a compact operator $\Tb=\Tb_{P,\AF,\Om}$ whose eigenvalues satisfy estimates similar to the ones in Section \ref{Sect.3}.

The conditions we impose on $\m,V$ are the following.

\begin{cond}\label{condition} If $2l<\Nb, $ $\mu$ satisfies \eqref{AR+} with $0<s<\Nb-2l;\,$ if
 $2l>\Nb,$ $\mu$ satisfies \eqref{AR-} with $0<s<\Nb$; the density $ V$ belongs to $L_{\vartheta,\m},$ \, $\vartheta=\max(1,\theta), \theta =\frac{s}{2l-\Nb+s}.$
\end{cond}

The  quadratic form is defined in  the following way. By results of Sect.2, under the above conditions,
the quadratic form $\sbb[v]=\int_\Om |v(X)|^2 P(dX),$ defined initially for $v\in H^l(\Om)\cap C(\Om),$ admits
 extension by continuity to the whole of $H^l(\Om)$ as a bounded quadratic form.
It its turn, operator $\AF$ maps $L_2(\Om)$ to $H^l(\Om),$ thus the composition
\begin{equation*}
    \tb[u]=\sbb[\AF u]
\end{equation*}
is well defined on $L_2(\Om)$ and determines a bounded operator $\Tb=\Tb_{P,\AF}.$

The action of operator $\Tb$ can be described explicitly, in a way similar to \cite{RIntegral}.
By the standard polarization, the sesquilinear form of operator $\Tb$ has the form
\begin{equation*}
    (\Tb u,v)_{L_2(\Om)}=\int (\AF u)(X)\overline{(\AF v)(X)}P(dX).
\end{equation*}
For a fixed $v\in L_2,$ the linear  functional $\psi_v(w)=\int {w(X)}\overline{{(\AF v(X))}}P(dX)$ is continuous on $H^l_0(\Om)$ and its norm is majorized
 by $||v||_{L_2}.$ This defines  the product  $\overline{(\AF v(X))}P$ as an element in the negative order Sobolev space of distributions $H^{-l}(\Om)$, and
 $v\mapsto \overline{(\AF v)} P$ is a continuous mapping from $L_2(\Om)$ to $H^{-l}(\Om).$ In its turn, for $h\in H^{-l}(\Om),$
the mapping  $u\mapsto \langle\AF u,h \rangle,$ defined initially on $u\in C_0^\infty(\Om),$ acts as $\vf_h:u\mapsto \langle u, \AF^* h\rangle,$ where
 $\AF^*:H^{-l}\to L_2(\Om)$ is the adjoint operator for $\AF:L_2(\Om)\to H^l_0,$ so $\AF^* h\in L^2(\Om)$ and the mapping $\vf_h$ extends by continuity to $u\in L^2(\Om)$.
  As a result, we obtain $(\Tb u,v)_{L_2(\Om)}=(\AF^*P\AF u,v)_{L_2(\Om)}$ and operator $\Tb$ factorizes as
  \begin{equation*}
    \Tb=\AF^* P \AF: L_2(\Om)\overset{\AF}{\longrightarrow}H^{l}_0(\Om)\overset{P}{\longrightarrow}H^{-l}(\Om)\overset{\AF^{*}}{\longrightarrow}L_2(\Om),
 \end{equation*}
all factors being continuous mappings between the corresponding spaces. This justifies our writing
$\Tb=\AF^*P \AF.$

Now we can formulate the eigenvalue estimates.
\begin{thm}\label{eigenv AF}Let $\Om$ be a bounded domain, $P=V \m$ be a singular measure in $\Om$, $\AF$ be a compactly supported pseudodifferential operator of order $-l$ and  Condition \ref{condition} be satisfied. Denote by $\Kb_\pm$ the quantity
\begin{gather}\label{AF coefficient}
    \Kb_{\pm}=\Ac(\m)^{\theta-1}\int_{\Om}V_{\pm}(X)^{\theta}\m(dX), \, {\mathrm{if}}\,\, 2l<\Nb,\\\nonumber
    \Kb_{\pm}=\Bc(\m)^{1-\theta^{-1}}\left(\int_{\Om}V_{\pm}(X)\m(dX)\right)^{\theta}\m(\Om)^{1-\theta}, \, {\mathrm{if}}\,\, 2l>\Nb.
\end{gather}
Then
\begin{equation}\label{AFest}
n_{\pm}(\la, \Tb)\le C(\AF,\Om)\Kb_{\pm}\la^{-\theta}, \la>0,
\end{equation}
and
\begin{equation}\label{AFsharperEst}
\limsup\la^{\theta}n_{\pm}(\la, \Tb)\le C(l,\Om)\sup_{(X,\X)\in \mathrm{S}^*\Om}|\ab_{-l}(X,\X)|^2\Kb_{\pm}.
\end{equation}
\end{thm}
  The rough estimate \eqref{AFest} in the first part of Theorem \ref{eigenv AF} is proved  quite similarly to Theorem 2.3 in \cite{RSh2} where the critical case was considered.
   Namely, denote by $\AF_0$ operator $(-\Delta_{\Dc}+1)^{-l/2}$ where $\Delta_{\Dc}$ is the Laplacian in a domain $\Om''\supset\Om'$ with Dirichlet boundary conditions; $\AF_0$ maps $L_2(\Om')$ to $H^{l}(\Om')$. Thus,
   \begin{equation*}
    \Tb_{P,\AF}=\AF^*P\AF= (\AF_0^{-1}\AF)^*(\AF_0^* P \AF_0)(\AF_0^{-1}\AF).
   \end{equation*}
  In this product, the middle term is exactly the operator $\Tb_{P,\AF_0}$ discussed  in Sect.3, for  which  the required eigenvalue estimates are already justified. Operator $\AF_0^{-1}\AF$ is bounded in $L_2(\Om),$ together with its adjoint, and this gives the first estimate. The second, sharper estimate \eqref{AFsharperEst} follows from the fact that the essential norm of $\AF_0^{-1}\AF$ equals $M(\AF):=\sup_{(X,\X)\in \mathrm{S}^*(\Om)}{|\ab_{-l}(X,\X)|}.$ Therefore, say, the estimate for positive eigenvalues follows from the operator inequality
  \begin{equation*}
    \Tb_{P,\AF}\le (M(\AF)+K)^*\Tb_{P,\AF_0}(M(\AF)+K),
  \end{equation*}
  with a compact operator $K.$ It is known that multiplication by the latter operator leads to a faster decay of eigenvalues, and this gives \eqref{AFsharperEst}.

   An even  sharper estimate is found later on, after some more localization properties are studied.

\subsection{Operators in $\R^\Nb$ and pseudodifferential CLR estimates} For operators in the whole Euclidean space we express the conditions for eigenvalue estimates in the terms of the mapping  properties of the pseudodifferential operator $\AF.$
\begin{thm}\label{ThmPsdoRN} Let $2l<\Nb$ and let  the order $-l$ pseudodifferential operator $\AF$ map $L_2(\R^\Nb)$ to $\Lb^l.$ Suppose that measure $\m$ satisfies \eqref{AR+} with some $s\in(\Nb-2l,\Nb)$ and $V$ is a real function, $V\in L_{\theta, \m},$ $\theta=\frac{s}{2l-\Nb+s}.$ Then for operator $\Tb=\Tb_{P,\AF},$ $P=V\m,$ defined by the quadratic form $\tb_{P,\AF}[u]=\int |(\AF u)(X)|P(dX)$ in $L_2(\R^\Nb)$ the estimate holds
\begin{equation}\label{EstPsDOE}
    n_{\pm}(\la, \Tb)\le C||\AF||_{L_2\to \Lb^l}^{2\theta}\Ac^{\theta-1}\int_{\R^{\Nb}}V_{\pm}(X)^{\theta}\m(dX) \la^{-\theta}.
\end{equation}
The estimate of the form \eqref{EstPsDOE} holds in the case of a complex-valued function $V$, with $V_\pm$ replaced by $|V|,$ for the counting function $n(\la, \Tb)$ of singular numbers of $\Tb.$
\end{thm}
\begin{proof} The estimate follows immediately from Theorem \ref{ThmRN} (see Remark \ref{RemRN}) due to the identity
\begin{equation*}
    \Tb_{P,\AF}=\AF^*P\AF=((-\D)^{l/2}\AF)^{*}\Tb_{P,(-\D)^{-l/2}}((-\D)^{l/2}\AF),
\end{equation*}
with $(-\D)^{l/2}\AF$ bounded in $L_2(\R^\Nb).$
\end{proof}

 From Theorem \ref{ThmPsdoRN} in the usual way follows the general CLR type estimate:
 \begin{cor}\label{cor.CLR}Let $\DF$ be an order $l<\Nb/2$ pseudodifferential operator in $\R^\Nb$ such that $\DF^{-1}$ maps $L_2(\R^\Nb) $ to $\Lb^{l}$. Let measure $\m$ satisfy \eqref{AR+} with $s\in(\Nb-2l,\Nb)$ and $V\in L_{\theta,\m},$ $\theta=\frac{s}{2l-\Nb+s}.$ Then for  the Schr\"odinger type operator $\HF=\DF^*\DF-P$  defined in $L_2(\R^\Nb) $ by means of the quadratic form $\hF[u]=\|\DF u\|_{L_2}^2-\int V|u(X)|^2\m(dX)$ the following estimate for the number $N_-(\HF)$ of negative eigenvalues holds
 \begin{equation*}
    N_-(\HF)\le C(\DF)\Ac^{\theta-1}\int_{\R^{\Nb}}V_+(X)^{\theta}\m(dX).
 \end{equation*}
 \end{cor}
Results for $2l>\Nb$ follow in the same way from Theorem \ref{Ex.BB}.
\subsection{Singular numbers estimates for non-selfadjoint operators.}\label{nonsa} We extend here the class of operators for which  spectral estimates are obtained. Let $\BFB$ be an order $-\vk$ pseudodifferential operator in  $\R^\Nb,$ $\vk<\Nb,$  and $\Mc$ be a closed set in $\R^\Nb.$ At this stage, we suppose  that $\Mc$ is compact. Let $\m$ be a Borel measure on $\Mc$ which satisfies condition \eqref{AR+} with some $s\in(\Nb-\vk, \Nb].$

We consider operators of the form
\begin{equation}\label{Tr1}
    \Gb=\Gb_{P_1,P_2,\BFB}=P_1\BFB P_2, \, P_j=V_j\m,
\end{equation}
in $L_{2,\m}(\Mc).$ Here $V_j$ are complex-valued $\m$- measurable functions on $\Mc,$ subject to conditions $V_j\in L_{r_j,\m}(\Mc)$ with some conditions imposed on $r_j,$ to be specified later.

Such operators have been systematically studied in \cite{ET}, Chapter 5, for the case $s=\Nb,$ i.e., for measure $\m$ absolutely continuous with respect to the Lebesgue measure, and later in \cite{Triebel2}, Chapter V, for an arbitrary $s,$  under the additional condition that $\mu$ coincides with the Hausdorff measure $\Hc^s $ of dimension $s$ on $\Mc, $  moreover the two-sided estimate \eqref{AR} was required.

The conditions imposed on $r_j$ in \cite{Triebel2} are the following (we present them in our notations, for the Hilbert space case $(p=2)$):
\begin{equation}\label{CondTr r}
    r_j>2; 1/r_1+1/r_2<\frac{\vk-\Nb+s}{s}=\theta^{-1}.
\end{equation}
Under these conditions, an estimate for the  eigenvalues $\la_k(\Gb), $ $|\la_1|\ge|\la_2|\ge\dots$ of operator $\Gb$ (counted with algebraic multiplicity), obtained in \cite{Triebel2}, see, e.g., Theorem 28.2 there,
sounds:
\begin{equation}\label{Est.Tr}
    |\la_k(\Gb)|\le C_{\ref{Est.Tr}}(\Mc,r_1,r_2)\|V_1\|_{L_{r_1,\m}}\|V_2\|_{L_{r_2,\m}}k^{-\theta^{-1}}, k\ge1.
\end{equation}
For the case of a self-adjoint operator $\Gb$, estimate \eqref{Est.Tr} is equivalent to a similar estimate, of the same order, for the singular numbers $s_k(\Gb).$ Both inequalities in \eqref{CondTr r} are not sharp: constant $C_{\ref{Est.Tr}}(\Mc,r_1,r_2)$ depends, due to \cite{Triebel2}, on the set $\Mc$ and deteriorates as $\Mc$ grows or  as $1/r_1+1/r_2$ approaches $\theta^{-1},$ or as one of $r_j$ approaches $2.$
In the case $1/r_1+1/r_2=\theta^{-1}<1$, called in \cite{ET}, \cite{Triebel2} the 'limiting' one, the conditions imposed upon one of the weight functions $V_j$ are strengthened, namely, say, for $V_1$, it is required that $V_1$ belongs to the Orlicz space, $L^{r_1}\log(1+L).$

We are going to demonstrate here that using our approach, these results, in the Hilbert space space setting, can be improved in several aspects. We show that the estimates for singular numbers, with constant \emph{not depending on}  $\Mc$ for the 'limiting' case, $1/r_1+1/r_2=\theta^{-1}$ follow from our general estimates, again, under the condition that only the upper estimate in \eqref{AR}, namely, \eqref{AR+} holds for measure $\m.$ In the discussion to follow, for simplicity of formulations, we consider the most important  case of $
\BFB=(-\Delta)^{-\vk/2}$ in $\R^\Nb$.
\begin{thm}\label{tr.improved} Let $\BFB=(-\Delta)^{-\vk/2}$ be the pseudodifferential operator in $\R^\Nb$, $\vk<\Nb$,  $\Mc$ be a compact set in $\R^\Nb$ supporting a measure satisfying \eqref{AR+}. Then $\Gb$ as an operator of the form \eqref{Tr1}, with $V_j\in L_{r_j,\m},$ $r_j>2,$ $1/r_1+1/r_2=\theta^{-1}=\frac{\vk-\Nb+s}{s},$ can be correctly defined as a bounded operator in $L_{2,\m}$ and for its singular numbers $s_k(\Gb)$ the estimate holds
\begin{equation}\label{Est.Tr.improved}
    s_k(\Gb)\le C\Ac^{1-\theta^{-1}}k^{-\theta^{-1}}\|V_1\|_{L_{r_1,\m}}\|V_2\|_{L_{r_2,\m}},
\end{equation}
where $\theta_j=\frac{s}{2l_j-\Nb+s}$, $l_j=\frac{\Nb-s}{2}+\frac{s}{r_j}<\Nb/2.$
\end{thm}
\begin{rem} In terms of the distribution function of singular numbers, \eqref{Est.Tr.improved} takes the form
\begin{equation}\label{Est.Tr.singular}
    n(\la,\Gb)\le C \la^{-\theta}\Ac^{\theta-1}\left[\int |V_1(X)|^{r_1}\m(dX)\right]^{\theta/{ r_1}}\left[\int |V_2(X)|^{r_2}\m(dX)\right]^{\theta/{r_2}}.
\end{equation}
\end{rem}
It is important to note that the constants in \eqref{Est.Tr.improved}, \eqref{Est.Tr.singular} do not depend on the set $\Mc.$ Therefore, similarly to Theorem \ref{ThmRN}, the above result extends to arbitrary measures satisfying \eqref{AR+}.
\begin{proof} We start by giving an exact definition of operator $\Gb.$ Since $\vk=2l_1+2l_2,$ consider the composition
\begin{equation}\label{Gb comp}
    \Gb=\VF_2 \VF_1^*,
\end{equation}
where $\VF_j=V_j \G_{\Mc}(-\Delta)^{-l_j/2}, $
 and $\G_{\Mc}$ is the operator of restriction of functions in the Sobolev space $H^{l_j}$ to $\Mc.$ Operator $\VF_j$ is considered as acting from $L_2(\R^{\Nb})$ to $L_{2,\m}(\Mc).$  As it follows from Proposition \ref{Sob.Embed} and the H\"older inequality, cf. Section \ref{2l<N.sub}, operator $\VF_j: L_2(\R^\Nb)\to L_{2,\m}(\Mc)$ defined initially on continuous functions admits a bounded extension by continuity to the whole of $L_2(\R^\Nb)$    under the conditions of our theorem. Additionally, similarly to Section 4.1, operator $\VF_j^*=(-\Delta)^{-l_j}\G_{\Mc}^*V_j: L_{2,\m}(\Mc)\to L_2(\R^\Nb)$ can be expressed by duality as the composition of continuous operators
 \begin{equation}\label{VF.factorize}
    \VF_j^*=L_{2,\m}(\Mc)\overset{V_j}{\longrightarrow}L_{\frac{2r_j}{2+r_j},\m}(\Mc)
    \overset{{\G_{\Mc}^*}}{\longrightarrow}\Lb^{-l_j}(\R^\Nb)\overset{(-\Delta)^{-l_j/2}}{\longrightarrow}L_2(\R^\Nb).
 \end{equation}
 Here $\G_{\Mc}^*:L_{\frac{2r_j}{2+r_j},\m}\to \Lb^{-l_j}$ is the embedding operator into the negative order homogeneous Sobolev space $\Lb^{-l_j}$ of distributions (dual to $\Lb^{l_j}$). This operator is continuous since it is adjoint to the embedding of $\Lb^{l_j}$ into $L_{\frac{2r_j}{r_j-2},\m},$ the latter being bounded by Proposition \ref{Sob.Embed}.
 Therefore, the operator $\Gb=\VF_2 \VF_1^*,$ defined in $L_{2,\m}$ by the sesquilinear form
 $\gb[u,v]= \langle\VF_1^* u,\VF_2^*v\rangle $  acts as the composition of continuous operators

 \begin{gather*}
    \Gb: L_{2,\m}(\Mc)\to L_{\frac{2r_2}{2+r_2},\m}(\Mc)
   \to\Lb^{-l_2}(\R^\Nb)\to L_2(\R^\Nb)\to\\\nonumber \to \Lb^{l_1}(\R^\Nb)\to L_{\frac{2r_1}{2-r_1},\m}(\Mc)\to L_{2,\m}(\Mc).
 \end{gather*}
 This reasoning  justifies the representation \eqref{Gb comp}.

 By the Ky Fan inequality, for the singular numbers of operators in \eqref{Gb comp}, $s_{2k-1}(\Gb)=s_{2k-1}(\VF_2 \VF_1^*)\le s_{k}(\VF_2)s_{k}(\VF_1^*)$
 For operator $\VF_2,$ we have $s_k(\VF_2)=s_k(\VF_2^*\VF_2)^{\frac12}=s_k(\Tb_2)^{\frac12},$
 where $\Tb_2=\Tb_{U_2\m,\AF_2}$ is the operator considered in Sect. 4, with $V=U_2=|V_2|^2$ and $\AF=\AF_2=(-\Delta)^{-l_j/2}.$ By Theorem \ref{ThmRN},
 \begin{equation*}
    n(\la, \Tb_2)\le C \Ac(\m)^{\theta_2-1}\int |U_2(X)|^{\theta_2} P(dX)\la^{-\theta_2}=
    C \Ac(\m)^{\theta_2-1}\int |V(X)|^{2\theta_2} P(dX)\la^{-\theta_2}.
 \end{equation*}
 The last estimate, expressed in terms of singular numbers, takes the form

\begin{equation}\label{est V1}
    s_k(\VF_2)=s_k(\Tb_2)^{\frac12}\le C \Ac(\m)^{(1-\theta_2^{-1})/2}\|V_2\|_{2\theta_2,\m} k^{-1/\theta_2}.
 \end{equation}
 Similarly, for $s_k(\VF_1^*)=s_k(\VF_1),$
 \begin{equation}\label{est V2}
  s_k(\VF_1)\le C \Ac(\m)^{(1-\theta_1^{-1})/2}\|V_1\|_{2\theta_1} k^{-1/\theta_2}.
 \end{equation}
 We combine \eqref{est V1} and \eqref{est V2} by means of the Ky Fan inequality to get
 \begin{equation*}
    s_{2k-1}(\Gb)\le s_k(\VF_1)s_k(\VF_2)\le C \Ac(\m)^{(2-\theta_1^{-1}-\theta_2^{-1})/2}\|V_1\|_{2\theta_1}\|V_2\|_{2\theta_2,\m}k^{-(\theta_1^{-1}+\theta_2^{-1})},
 \end{equation*}
 which, since $\theta_1^{-1}+\theta_2^{-1}=\theta^{-1}$, coincides with \eqref{Est.Tr.improved}.
\end{proof}
We skip the discussion of  spectral estimates for other relations between $r_1,r_2,\vk,\Nb,$ which are established by means of a similar reasoning, using the results of  our Section 4.

Another type of non-self-adjoint operators with singular weights, considered in \cite{Triebel2}, are
\begin{equation*}
    \Gb={\AF_1 P \AF_2},
\end{equation*}
where $\AF_1,\AF_2$ are pseudodifferential operators in $\R^\Nb$ of orders $-l_1,-l_2$ and $P$ is a singular measure of the form $P=V\m,$ with $\m$ being, as before, the dimension $s$ Hausdorff measure on a set $\Mc\subset \R^\Nb$, satisfying \eqref{AR}, and $V$ being a $\m$-measurable function on $\Mc$ belonging to $L_{r,\m}.$ As in the previous case, the conditions imposed in \cite{Triebel2} on the parameters of this operator are not sharp, which leads to non-sharpness of spectral estimates.

Spectral (singular numbers) estimates for such operators are obtained here in the same way as in Theorem \ref{tr.improved}. Namely we perform the factorization, this time of the weight function $V$, which leads to the factorization of the operator $\Gb.$ Again, we present the most aesthetic case.

\begin{thm}\label{Th.conversed} Let $\AF_j=(-\Delta)^{l_j/2},\, j=1,2,$ in $\R^\Nb$, $2l_j<\Nb$, $2l=l_1+l_2.$  Let, further, $\m$ be a singular measure satisfying \eqref{AR+} with exponent $s,$ $\Nb-2l_j<s\le \Nb,$ and $P=V\m,$ $V\in L_{\theta,\m}.$ We set $\theta_j=\frac{s}{2l_j-\Nb+s}>1,$ so that  $2\theta^{-1}=\theta_1^{-1}+\theta_2^{-1}= \frac{s}{2l-\Nb+s},$ $\theta>1.$
Consider operator $\Gb=\Gb_{j_1,j_2,P}=\AF_2 P\AF_1$ in $L_2(\R^\Nb)$ defined by the sesquilinear form
\begin{equation*}
    \gb[u,v]=\int_{\R_\Nb}(\AF_1 u)(X)\overline{(\AF_2 v)(X)}P(dX), \, u,v\in L_2(\R^\Nb.)
\end{equation*}
Then for the singular numbers of $\Gb$ the estimate holds
\begin{equation}\label{Eq.Th.conversed}
    n(\la, \Gb)\le C \int{|V(X)|^{\theta}}\m(dX)\la^{-\theta.}
\end{equation}
\end{thm}
\begin{proof} As before, it suffices to consider the case of real, non-negative $V.$ We represent $V$ as $V=V_1V_2,$ $V_j=V^{\g_j},$ where $\g_j=\theta_{2-j}(\theta_1+\theta_2)^{-1}=\theta/(2\theta_j),$ $\theta_j=\frac{s}{2l_j-\Nb+s},$ so that $\g_1+\g_2=1,$ $\theta_1^{-1}+\theta_2^{-1}=(2\theta)^{-1}.$
With this representation, operator $\Gb$ factorizes as
\begin{equation*}
    \Gb=\UF_2^*\UF_1, \, \UF_j=V_j\G \AF_j: L_2(\R^\Nb)\to L_{2,\m}(\Mc).
\end{equation*}
As in Theorem \ref{tr.improved}, it follows from Lemma \ref{lem.RN} that operators $\UF_j$ are bounded as acting from
$L_2(\R^\Nb)$ to $L_{2,\m}(\Mc)$ and are factorized as compositions of continuous operators.
By the Ky Fan inequality
\begin{equation}\label{KF.tr}
    n(\la_1\la_2,\UF_2^*\UF_1 )\le n(\la_1,\UF_1)+n(\la_2,\UF_2),
\end{equation}
where $\la_j,$ $\la_1\la_2=\la,$ will be fixed later. For the  terms on the right in \eqref{KF.tr},  we have
\begin{equation*}
    n(\la_j,\UF_j)=n(\la_j^2,\UF_j^*\UF_j).
\end{equation*}
Operator $\UF_j^*\UF_j$ coincides with $\AF_j^* V_j^2 \AF_j,$ and by Theorem \ref{ThmRN} the estimate holds
\begin{equation}\label{Tr.Uj}
    n(\la_j,\UF_j)\le C \Ac^{\theta_j-1}\la^{-2\theta_j}\int V_j(X)^{2\theta_j}\m(dX)=C \Ac^{\theta_j-1}\la^{-2\theta_j}\int V(X)^{\theta}\mu(dX),
\end{equation}
since $V_j^{2\theta_j}=V^{2\g_j\theta_j}=V^{\theta}.$
Finally, we set $\la_j=\la^{\frac{\theta_j}{\theta_1+\theta_2}},$ so that $\la_1\la_2=\la,$
and now \eqref{Tr.Uj}, \eqref{KF.tr} give the required inequality \eqref{Eq.Th.conversed}.
\end{proof}
\begin{rem}In a similar way, other results of Sections 3,4 are carried over to non-selfadjoint operators.
\end{rem}

\section{Localization} In this section we present some auxiliary  results on perturbations and localization in eigenvalue estimates  for operators of the form $\Tb_{P,\AF}$, needed for obtaining eigenvalue asymptotics. Some of these results are analogous to the ones in Sect.3 in \cite{RSh2}. We, however, give complete proofs.

It is convenient to describe  perturbation and localization properties without supposing that asymptotic formulas hold. So, the asymptotic characteristics of eigenvalue distribution, the ones introduced in \cite{BS}, are used.
\begin{defin}\label{As.char}Let $\Tb$ be a compact self-adjoint operator. For $\theta>0,$ we define the quantities
\begin{equation}\label{DefinAsCharD}
    \Db_{\pm}^{\theta}(\Tb)=\limsup_{\la\to 0}n_{\pm}(\la,\Tb)\la^{\theta},\, \, \db_{\pm}^{\theta}(\Tb)=\liminf_{\la\to 0}n_{\pm}(\la,\Tb)\la^{\theta},
\end{equation}
where, recall, $n_{\pm}(\la,\Tb)$ is the distribution function of positive (negative) eigenvalues of $\Tb.$
\end{defin}
Of course, these quantities equal zero or infinity for all values $\theta$ except, possibly, just one. However, if for some $\theta,$ for a certain choice of $\pm$ sign, both $ \Db_{\pm}^{\theta}(\Tb)$ and $\db_{\pm}^{\theta}(\Tb)$ are finite, nonzero, and equal, this means that the eigenvalues of $\Tb$ of the corresponding sign are subject to an asymptotic formula of order $\theta$.
Some of the   statements to follow appear mostly  in two versions each, the subcritical and supercritical ones. The proofs are identical, up to notations, therefore we present them only for one statement in the pair.

We will systematically refer to the following elementary observation.
\begin{proposition}\label{order monotonicity} Let $l<l'$ and pseudodifferential operator $\AF'$ have order $-l'$; for $2l<\Nb$ it is supposed that $2l'<\Nb$. Suppose that measure $\m$ and density $V$ satisfy Condition \eqref{condition}. Then $n_{\pm}(\la,\Tb_{V\m,\AF'})=o(\la^{-\theta})$, $\theta=\frac{s}{2l-\Nb+s}$.
\end{proposition}
This property follows from the fact that for  fixed $P, s,$ the exponent $\theta$ grows as the order $l$ grows.

\begin{lem}[Lower order perturbation]\label{perturbation}
Let $\AF$ be a pseudodifferential operator in $\Om$ of order $-l$ and $\BFB$ be a pseudodifferential operator of order $-l'<-l$, $\AF_1=\AF+\BFB$. Let Condition \ref{condition} be satisfied. Then for $\theta=\frac{s}{s+2l-\Nb},$
\begin{equation}\label{Pert.loworder}
    \Db_{\pm}^\theta(\Tb_{P,\AF})=\Db_{\pm}^{\theta}(\Tb_{P,\AF_1}), \, \db_{\pm}^\theta(\Tb_{P,\AF})=\db_{\pm}^{\theta}(\Tb_{P,\AF_1}).
\end{equation}
In particular, if for the positive (negative) eigenvalues of operator $\Tb_{P,\AF}$ an asymptotic formula of order $\theta$ is valid, it is valid for eigenvalues of operator $\Tb_{P,\AF_1}$ as well, with the same coefficient.
\end{lem}
\begin{proof} Consider, e.g.,  the subcritical  case. The difference $\Tb_{P,\AF}-\Tb_{P,\AF_1}$ is defined by the quadratic form which is the sum of terms $2\int \re\left[(\BFB u)\overline{(\AF u)}\right]P(dX),$    and $\int|\BFB u(X)|^2P(dX).$ For the latter form,  $\BFB$ is a lower order perturbation and Proposition \ref{order monotonicity} applies. As for the first form above, we have
\begin{gather}\label{loworder2}
    \left| \int \re\left[(\BFB u)\overline{(\AF u)}\right]P(dX)\right|\le 2\e\int|\AF u(X)|^2|V(X)|\m(dX)\\\nonumber +2\e^{-1}\int|\BFB u(X)|^2|V(X)|\m(dX).
\end{gather}
Here, the first term generates operator $2\e \Tb_{|V|\m, \AF}$ with an arbitrarily small $\e$, and the second term generates the operator already considered, with eigenvalues subject to $n_{\pm}(\la)=o(\la^{-\theta})$;  by the arbitrariness of $\e$,  $ n_{\pm}(\la,\Tb_{P,\AF}-\Tb_{P,\BFB})=o(\la^{-\theta})$.
\end{proof}

\begin{lem}[Cf. Observation 4 in \cite{RSh2}]\label{Low.Order.Lem}
 Suppose that  operator $\AF$ is of order $-l$ and measure  $P=V\m$ satisfies Condition \eqref{condition}. Let $\hi$ be the characteristic function of a closed set $\Eb\subset\Om$ such that $\Eb\cap\Mc=\varnothing.$   Then for  operator $\Tb_{P,\AF \hi}$ defined by the quadratic form
\begin{equation}\label{Loc2}
    \int_{\Om}|\AF(\hi u)(X)|^2 P(dX),
\end{equation}
the eigenvalues satisfy
\begin{equation}\label{Loc3}
    \Db_{\pm}^{\theta}(\Tb_{P,\AF \hi})=0.
\end{equation}
\end{lem}
This result demonstrates the following spectral localization property: if the quadratic form is restricted to functions supported away from the support of the measure $\m$, the eigenvalues decay faster than this might be prescribed by the general estimate.
\begin{proof}Let $\psi\in C_0^\infty(\Om)$ be a function that is equal $0$ on $\Eb$ and $1$ in a neighborhood of $\Mc$. Then
 \begin{equation*}
 \tb_{P,\AF\hi}[u]=\int V(X)|\psi(X)\AF(u \hi )(X)|^2\m(dX)=\int V(X)|([\psi,\AF](u\hi))(X)|^2\m(dX),
\end{equation*}
and since the commutator  $[\psi,\AF]$ has order $-l-1,$ Lemma \ref{Low.Order.Lem} applies.
\end{proof}

It follows, in particular, that the eigenvalue counting function gets a lower order perturbation if  operator $\AF$ is perturbed outside a neighborhood of the set $\Mc$. This circumstance gives us freedom in choosing cut-off functions away from the  $\Mc$ or adding operators that are smoothing near $\Mc$ -- the possibility already mentioned.

A complication which is often encountered in the study of eigenvalue asymptotics is the non-additivity of asymptotic coefficients: if for operators $T_1$, $T_2$ asymptotic formulas of the same order for $n(\la)$ are known, this does not automatically imply a similar formula for $T_1+T_2$ (and does not generally imply any asymptotic formula for eigenvalues at all.) In the following  statement,  important in the study of eigenvalue asymptotics, it is shown that if two measures have supports separated by a positive distance, then, up to a lower order term, the counting functions behave additively with respect to the measures: measures are, so to say,  spectrally almost orthogonal.

\begin{lem}\label{LemLocalization}
Let Condition \ref{condition} be satisfied, $\AF$ be an operator of order $-l$.
 Suppose that  $P=P_1+P_2$, where $P_j=V_j \m_j$ is a measure supported on a compact set $\Mc_j$ and  $\dist(\Mc_1, \Mc_2)>0$. Then
\begin{equation}\label{Localization}
    n_{\pm}(\la,\Tb_{P_1+P_2})=n_{\pm}(\la, \Tb_{P_1})+n_{\pm}(\la, \Tb_{P_2})+o(\la^{-\theta})\ \mbox{ as }\ \la\to 0,
\end{equation}
in particular,
\begin{equation}\label{LocDd}
    \Db^\theta_{\pm}(\Tb_{P_1+P_2})\le \Db^\theta_{\pm}(\Tb_{P_1})+\Db^\theta_{\pm}(\Tb_{P_2}),\,
    \db^\theta_{\pm}(\Tb_{P_1+P_2})\ge \db^\theta_{\pm}(\Tb_{P_1})+\db^\theta_{\pm}(\Tb_{P_2})
\end{equation}
\end{lem}
\begin{proof} Consider two disjoint open sets $\Om_1,\Om_2\subset \Om$, such that $\overline{\Mc_j}\subset\Om_j$ and $\overline{\Om_1}\cup\overline{\Om_2}\supset\Om$.
Every function $u\in L_2(\Om)$ splits into the (orthogonal) sum $u=u_1\oplus u_2$, $u_j\in L_2(\Om_j)$. The quadratic form of the operator $\Tb_{P_1+P_2}$ splits as follows
\begin{gather}\label{splitting}
  \tb_{P_1+P_2}[u] :=  (\Tb_{P_1+P_2}u,u)_{L^2(\Om)}  =\int_\Om V_1(X)|\AF(u_1)(X)+ \AF(u_2)(X)|^2 \m_{1}(dX)\\\nonumber
+  \int_\Om V_2(X)|\AF(u_1)(X)+ \AF(u_2)(X)|^2\m_2(dX) = \tb_{P_1}[u_1]+\tb_{P_2}[u_2]+\tb_R[u_1,u_2]:=\\
\nonumber  \int_\Om V_1(X)|\AF(u_1)(X)|^2 \mu_{1}(dX)+\int_\Om V_2(X)|\AF(u_2)(X)|^2\mu_{2}(dX)+\tb_R[u_1,u_2].
\end{gather}
The remainder term $ \tb_R[u_1,u_2]$ is a quadratic form in function $u=u_1\oplus u_2,$ $u_j\in L_2(\Om_j)$ with the following property: if a term in $\tb_R$ contains measure $P_j$, then it necessarily contains function
$u_{3-j}$, so it always contains a measure and a function with disjoint supports. If such a term has the form $\int_{\Om} V_1|\AF u_2|^2 \mu_1(dX)$, the corresponding operator $\Tb$ satisfies $n_{\pm}(\la, \Tb)=o(\la^{-\theta})$ by Lemma \ref{Low.Order.Lem}. If, on the other hand, such term has the form $\int_{\Om} V_1(\AF u_1)\overline{(\AF u_2)}\m_{1}(dX)$, then by the Schwartz inequality,
 $$\left|\int_{\Om} V_1(\AF u_1)\overline{(\AF u_2)}\m_1(dX)\right|\le \left(\int_{\Om} |V_1||\AF u_1|^2{\mu_{1}(dX)}\right)^{1/2}\left(\int_{\Om} |V_1||\AF u_2|^2\mu_{1}(dX)\right)^{1/2},$$
 and the last factor, again by Lemma \ref{Low.Order.Lem}, provides the required lower order estimate for the eigenvalue counting function.
Now we observe that  forms $\tb_{P_1}[u_1],\,\tb_{P_2}[u_2]$ in \eqref{splitting} act on orthogonal subspaces $L_2(\Om_1), L_2(\Om_2)$. Therefore the spectrum of the sum of the corresponding operators $\Tb_1,\Tb_2$  equals the union of the spectra of the summands, and hence $n_{\pm}(\la,\Tb_1+\Tb_2)=
n_{\pm}(\la,\Tb_1)+n_{\pm}(\la,\Tb_2).$ The term $\tb_R$ in  \eqref{splitting} makes a weaker contribution,
\begin{equation}\label{splitting1}
n_{\pm}(\la,\Tb_{P_1+P_2})=
n_{\pm}(\la,\Tb_1)+n_{\pm}(\la,\Tb_2)+o(\la^{-\theta}).
\end{equation}
Now consider  operator $\Tb_{P_1}$. It is defined by the quadratic form
 \begin{equation*}
 \tb_{P_1}[u]=\int V_1(X)|\AF(u_1\oplus u_2)(X)|\mu_{1}(dX).
 \end{equation*}

  Similarly to \eqref{splitting}, we represent it as
\begin{equation}\label{splitting2}\tb_{P_1}[u]=\int_\Om V_1(X)|\AF(u_1)(X)|^2 \mu_{1}(dX)+\tb_{R_1}[u_1,u_2],
\end{equation}
with $\tb_{R_1}$ having the same structure as $\tb_{R}$ in \eqref{splitting}. Again, form $\tb_{R_1}$ generates an operator  $\tb_{R_1}$ with eigenvalues satisfying $n_{\pm}(\la,\Tb_{R_1})=o(\la^{-\theta}),$
and we obtain
\begin{equation}\label{splitting3}n_{\pm}(\la, \Tb_{P_1})=n_{\pm}(\la,\Tb_1)+o(\la^{-\theta}).
\end{equation}
In the same way,
\begin{equation}\label{splitting4}n_{\pm}(\la, \Tb_{P_2})=n_{\pm}(\la,\Tb_2)+o(\la^{-\theta}).
\end{equation}
Finally, we substitute \eqref{splitting3}, \eqref{splitting4} into \eqref{splitting1} to obtain \eqref{Localization}, and therefore, \eqref{LocDd}. (Different signs in two inequalities in \eqref{LocDd} arise, of course, due to the fact that $\limsup(f(\la)+g(\la))\le \limsup f(\la)+\limsup g(\la),$ while $\liminf(f(\la)+g(\la))\ge \liminf f(\la)+\liminf g(\la).$ )
\end{proof}

An important corollary of this lemma provides us with the additivity property for asymptotic formulas.
\begin{cor}\label{Cor.addit} Under the conditions of Lemma \ref{LemLocalization}, if, for some sign, asymptotic formulas
for eigenvalues
\begin{equation}\label{As.j}\lim_{\la\to 0}\la^\theta n_{\pm}(\la, \Tb_{P_j})=A_j^{\pm}, \, j=1,2,
\end{equation}
hold, then a similar formula is valid for $\Tb_{P_1+P_2}=\Tb_{P_1}+\Tb_{P_2}:$
\begin{equation}\label{As.1+2}\lim_{\la\to 0}\la^\theta n_{\pm}(\la, \Tb_{P_1}+\Tb_{P_2})=A_1^{\pm}+A_2^{\pm}.
\end{equation}
\end{cor}
Relations \eqref{As.j} are equivalent to $\Db_\pm^\theta(\Tb_{P_j})=\db_\pm^\theta(\Tb_{P_j})=A_j^\pm$,
and now \eqref{As.1+2} follows from \eqref{LocDd}.

Another corollary of Lemma \ref{LemLocalization} allows us to separate the positive and  negative parts of the function $V$ when studying the distribution of the positive and  negative eigenvalues of $\Tb_{P} $ separately.
\begin{cor}\label{cor.plusminus}
 Suppose that Condition \ref{condition} is satisfied.
 Let $\Mc$ be a compact set, $\m=\m_++\m_-,$ $\supp\mu_\pm=\Mc_{\pm},$ $\dist(\Mc_+,\Mc_-)>0$, and $V\in L_{\vartheta}(\Mc).$ Let
$V_{\pm}\ge 0$ in $\Mc_{\pm}$, $V=0$ in $\Mc\setminus(\Mc_+\cup\Mc_-),$  $P=V\mu$, $P_{\pm}=V_{\pm}\m_{\pm}$.
Then
\begin{equation}\label{sign.separation}
n_{\pm}(\la, \Tb_P)=n_+(\la, \Tb_{V_\pm\mu})+o(\la^{-\theta}) \ \mbox{ as }\ \la\to 0,
\end{equation}
in particular, $\Db_\pm^\theta(\Tb_{P})=\Db_\pm^\theta(\Tb_{P_{\pm}}),$ $\db_\pm^\theta(\Tb_{P})=\db_\pm^\theta(\Tb_{P_{\pm}}),$
\end{cor}
In other words, up to a lower order remainder, the asymptotic behavior of the positive, resp., negative, eigenvalues of  operator $\Tb_{V\m}$ is determined by the positive, resp., negative, part of   density $V,$ as soon as these parts are separated. To prove this property, we can use \eqref{splitting1}, taking as $P_1$ the restriction of  measure $P$ to the set $\Mc_+$, and as $P_2$ its restriction to $\Mc_-$, and recall that $n_-(\la, \Tb_{V\mu,\Om_+})=n_+(\la, \Tb_{V\m,\Om_-})=0,$ where $\Om_{\pm}$ are neighborhoods of $\Mc_{\pm}$ such that $\dist(\Om_+,\Om_-)>0 $.

Finally, we can get rid of the condition for these sets to be separated.
\begin{thm}\label{PMseparation} Let measure $\m$, density $V$ and operator $\AF$ satisfy Condition \ref{condition}. Then \eqref{sign.separation} holds.
\end{thm}
\begin{proof} We follow mostly the reasoning in \cite{RSh2}, where a similar property was established in a more restricted form.
For a given $\ve>0$ we approximate  density $V$ by a  function $V_\ve$, continuous on $\Mc$ so that  $\int_{\Mc}|V-V_\ve|^\vartheta  \m(dX)<\ve$, $\vartheta=\max(\theta,1).$
By estimates in Sect. \ref{Sect.3}, $\Db^{\pm}_{\theta}(\Tb_{V}-\Tb_{V_{\ve}})\le C \ve^{\vartheta/\theta} $. Consider the set $\Mc_0=\{X\in \Mc: V_\ve=0\}$ and its $\de$ -neighborhood $\Mc_\de$. For sufficiently small $\de,$ the quantity $\int_{\Mc_{\de}}|V_{\ve}|^{\vartheta}\m(dX)$ is less than $\ve$. So, for the weight $V_{\ve}',$  the restriction of $V_{\ve}'$ to $\Mc\setminus\Mc_\de,$ by the same estimates in Sect \ref{Sect.3}, $\Db^{\pm}_{\theta}(\Tb_{V}-\Tb_{V'_{\ve}})\le C \ve^{\vartheta/\theta}.$ Therefore,
\begin{equation}\label{separation}
    \Db_\theta^{\pm}(\Tb_{V}-\Tb_{V'_{\ve}})< C\ve^{\vartheta/\theta}.
\end{equation}
Now, density $V'_{\ve}$ satisfies conditions of Corollary \ref{cor.plusminus}, and \eqref{sign.separation} follows from \eqref{separation} and the arbitrariness of $\ve$ by means of the Ky Fan inequality.
\end{proof}

\section{Eigenvalue asymptotics}\label{Sect.Asymp.}Derivation of eigenvalue asymptotic formulas follows now the pattern of \cite{RSh2} and \cite{RIntegral}.

\subsection{Asymptotic perturbation lemma}\label{SubS.BSLemma}For convenience of the Reader, we reproduce here the asymptotic perturbation lemma by M. Birman and M. Solomyak, see, e.g. \cite{BS}, Lemma 1.5, which we use systematically.
\begin{lem}\label{BSLemma} Let $T$ be a compact self-adjoint operator. Suppose that for any $\ve>0$ operator $T$ can be represented as the sum, $T=T_\ve+T_{\ve'}$, so that for $T_\ve$ the asymptotic formula for positive, resp., negative, eigenvalues holds, $\lim_{\la\to 0}\la^{\theta}n_{\pm}(\la,T_\ve)=\Cb^{\pm}_{\ve}$, for and for $T_{\ve}'$ both  relations
\begin{equation}\label{SmallPert}
    \lim_{\ve\to 0} \Db_\theta^{\pm}(T'_{\ve})=0
\end{equation}
hold. Then the limit $\Cb^{\pm}=\lim_{\ve\to 0}\Cb^{\pm}_{\ve}$ exists and for operator $T$ the asymptotic formula is valid, $\lim_{\la\to 0}\la^\theta n_{\pm}(\la,T)=\Cb^{\pm}.$
\end{lem}

\subsection{Compact Lipschitz surfaces}\label{subsect.Lip} The first result concerns measures supported on Lipschitz surfaces. A compact Lipschitz surface is, locally,
 a graph of a Lipschitz vector-function.

Let $\Si\subset\R^\Nb$ be a compact Lipschitz surface   of dimension $d: 0<d<\Nb$ and codimension $\dF=\Nb-d,$ in a domain $\Om\subset\R^{\Nb},$ defined, locally, in proper local
 co-ordinates $X=(\xb,\yb), \xb\in \R^d, \yb\in \R^{\dF},$  by the equation $\yb=\pmb{\f}(\xb)$, with Lipschitz vector-function $\pmb{\f}$: $|\pmb{\f}(\xb)-\pmb{\f}(\xb')|\le C_{\pmb{\f}}|\xb-\xb'|,$ $\xb\in\omega\subset\R^d$. Measure $\m$ on $\Si,$ generated by the embedding of $\Si$ into $\R^{\Nb}$ coincides with the $d$-dimensional Hausdorff
  measure $\Hc^d$ on $\Si$. For a compact Lipschitz surface, in one co-ordinate neighborhood, and therefore, globally, Ahlfors condition \eqref{AR} is satisfied, with $s=d.$ For  noncompact Lipschitz surfaces, the situation is somewhat more complicated, see Section \ref{Sect.noncompact}.

    By the Rademacher theorem, for $\m$-almost every  $X\in\Si$, there exists a tangent space $\Tt_X\Si$ to $\Si$ at the point $X$ and, correspondingly, the normal space $\Nt_X\Si$ which are identified naturally with the cotangent and the conormal spaces. By $\St_X\Si$ we denote the sphere $|\x|=1$ in $\Tt_X\Si.$
 \begin{thm}\label{Th.As.RSh} Let  real  function  $V$ on $\Si$ belong to $L_{\vartheta,\m}(\Si)$,
 $\vartheta=\max(\theta, 1),$ $\theta =\frac{d}{2l-\Nb+d}=\frac{d}{2l-\dF}$, $\dF<2l$ and $P=V\m$. Let $\AF$ be an order
  $-l$ pseudodifferential operator  compactly supported in $\Om\subset\R^\Nb$, with principal symbol $\ab_{-l}(X,\X)$.
At those points $X\in \Si$,  where the tangent plane exists, we define the auxiliary symbol $\rb_{-m}(X,\x)$, $\x\in \mathrm{T}_X\Si$, of order $-m=\dF-2l<0$,
\begin{equation}\label{r_-q}
    \rb_{-m}(X,\x)=(2\pi)^{-\dF}\int_{\mathrm{N}_X\Si}|\ab_{-l}(X,\x,\y)|^2 d\y, (X,\x)\in \mathrm{T}^*\Si,
\end{equation}
and the density
\begin{equation}\label{ro} \ro_{\AF}(X)=\int_{\mathrm{S}_X\Si} \rb_{-m}(X,\x)^{\theta}d\x.
\end{equation}

Then for the eigenvalues of  operator $\Tb_{V,\m,\AF}=\AF^{*}P\AF,$ $P= V\m$, the asymptotic formulas are valid
\begin{equation}\label{AsLambda}
     n_{\pm}(\la,\Tb_{V\m,\AF})\sim \la^{-\theta} A_{\pm}(V,\Si,\AF), \, \la\to 0,
    \end{equation}
where
  \begin{gather}\label{coeff}
   A_{\pm}(V,\Si,\AF)= \frac{1}{d (2\pi)^{d-1}}\int_\Si\int_{\St_X\Si} V_{\pm}(X)^\theta\rb_{-m}(X,\x)^{\theta}d\x \m(dX)= \\\nonumber
   \frac{1}{d (2\pi)^{d-1}} \int_\Si V_{\pm}(X)^{\theta}\ro_{\AF}(X) \m(dX), \end{gather}
with  density $\ro_{\AF}(X)$ defined in \eqref{ro}.
\end{thm}

In the particular case of  $\AF=\AF_0=(1-\D)^{-l/2}$ in a neighborhood of $\Si$ in $\R^\Nb$, we have  $\ab_{-l}(X,\X)=|\X|^{-l}$ and

\begin{gather*}
\rb_{-m}(X,\x)= (2\pi)^{-\dF}\int_{\R^{\dF}}(|\x|^2+|\y|^2)^{-l}d\y =\\
\nonumber |\x|^{-m}(2\pi)^{-\dF}\pmb{\om}_{\dF-1}\int_{0}^{\infty}\z^{\dF-1}(1+\z^2)^{-l}d\z =
\pmb{\om}_{\dF-1}\frac{1}{2(2\pi)^\dF}\Bb\left( \frac{\dF}{2}, \frac{l-\dF}{2}\right) |\x|^{-m},
\end{gather*}
where $\pmb{\om}_{\dF-1}$ is the volume of the unit sphere in $\R^{\dF}$, $\Bb$ is  Euler's  Beta-function.
Therefore,
\begin{equation}\label{ro.special}
\ro_{\AF_{0}}\equiv\pmb{\ro}(\dF,l)=\left[\frac{\om_{\dF-1}}{2(2\pi)^\dF}\Bb\left( \frac{\dF}{2}, \frac{l-\dF}{2}\right)\right]^{\theta}
\end{equation}
and, finally,
\begin{equation}\label{integral}
   n_{\pm}(\la, \Tb_{V,\AF_0})\sim \la^{-\theta}\Zb(d,\dF)\int_{\Si} V_{\pm}^{\theta}(X) \m(dX),
\end{equation}
\begin{equation}\label{IntCoeff}
\Zb(d,\dF)=\frac{\pmb{\om}_{d-1}}{d}\pmb{\ro}(\dF,l).
\end{equation}

\begin{proof}{The {Proof}}
 of  Theorem \ref{Th.As.RSh} follows closely the one of  Theorem 2.4 in \cite{RSh2}, with modifications  caused by  a different order of the operators involved. We give detailed  explanations of  the corresponding changes, directing  interested Readers to \cite{RSh2} for more details.

First, by the usual localization (see, e.g. how this procedure is performed in \cite{BSSib} or \cite{RT1}), it is sufficient to prove the asymptotic formula for $V$ supported in just one co-ordinate neighborhood of $\Si.$ In such case,
we can  approximate $V$ by a  density $V_{\ve}$, defined and smooth in a neighborhood of $\Si$ (therefore, in $\R^{\Nb}$) such the $L_{\vartheta,\m}$- norm of $V_{\ve}-V$ is small, less than $\ve^{\vartheta/\theta}$. By eigenvalue estimates in Sect. \ref{Sect.4}, the eigenvalue distribution functions for operators $\Tb_{V,\m,\AF}$ and $\Tb_{V_{\ve},\m,\AF}$ differ asymptotically  by less than $C \ve\la^{-\theta}$.  By the  asymptotic perturbation  Lemma \ref{BSLemma}, such approximation enables us to restrict our task to proving asymptotic formulas for such nice densities $V_{\ve}$ only,  passing then to limit as $V_{\ve}$ approaches $V$ in the $L_{\vartheta,\m}$ norm.
 On the next step, we separate the positive and negative eigenvalues of our operator, using Theorem \ref{PMseparation}.  In this way, the problem is reduced to
the case of a non-negative density  $V_{\ve}$, which, again after adding a perturbation with arbitrarily small $L_{\theta,\m}$-norm and using Lemma \ref{BSLemma}, we may suppose being the restriction to $\Si$ of a smooth compactly supported non-negative function in $\R^\Nb$; we drop the subscript $\ve$ further on. So,
$V=U^2,$ $U\in C_0^{\infty}(\R^\Nb)$, see \cite{RSh2} for the detailed description of this construction.

  Next, the spectral problem is reduced to the study of eigenvalues of an integral operator on $\Si$ with kernel having a weak singularity at the diagonal. This is done in the following way. Operator $\Tb_{V\m,\AF}$ is defined by the quadratic form $\tb_{V\m,\AF}[u],$ see \eqref{Pseu.1},
  which can be re-written  as
   \begin{equation*}
    (\Tb_{V\m,\AF} u,u)_{L_2(\Om)}=\tb_{V\m,\AF}[u]=(\G_\Si U\AF u,\G_\Si U\AF u)_{L_{2,\m}}=
    ((\G U\AF)^*(\G U\AF)u,u)_{L_2(\Om)},
   \end{equation*}
   where $\G:H^l(\Om)\to L_2(\mu)$ is the operator of restriction to $\Si.$
  In this way, operator $\Tb_{V\m,\AF}$ factorizes as
\begin{equation}\label{factorizationSi}
 \Tb_{V_{\ve}\m,\AF}=\KF^*\KF, \, \KF=\G_\Si U\AF: L_2(\Om)\to L_2(\Si,\m), \KF^{*}:L_2(\Si,\m)\to L_2(\Om).
 \end{equation}
 We know, however, that the nonzero eigenvalues of  operator $\KF^*\KF$ coincide with nonzero eigenvalues of $\RF=\KF \KF^*,$ counting multiplicities. By \eqref{factorizationSi},  operator $\KF \KF^*$ acts in $L^2(\Si,\m)$ as

\begin{equation}\label{transposed}
\RF=\KF \KF^*=   \G_{\Si}U\AF\AF^*U\G_{\Si}^*.
\end{equation}

Operator $U\AF\AF^*U$ is an order $-2l$ pseudodifferential operator in $\Om$ with principal symbol $\Rc_{-2l}(X,\X)=V(X)|a_{-l}(X,\X)|^2$, or, equivalently, it is a self-adjoint integral operator with kernel $R(X,Y,X-Y)$, smooth for $X\ne Y$ and having a weak singularity as $X-Y\to 0$. The leading singularity of this kernel, being the Fourier transform of the symbol  $\Rc_{-2l}(X,\X)$ in $\X$ variable, has the following structure. If $m=2l-\Nb$ is \emph{not} a positive even integer, the  leading singularity in $X-Y$ of
 $R(X,Y,X-Y)$ has the form $R_m(X,X-Y),$ smooth for $X\ne Y$ and  positively homogeneous of order $m$ in $X-Y$. In case when $m=2l-\Nb$ \emph{is} an even positive integer, then, in addition to the above term, there may be present a term with leading singularity of the form $R_{\log}(X,X-Y)\log|X-Y|,$ where $R_{\log}(X,X-Y)$ is a smooth function of all variables, being a homogeneous polynomial in $X-Y$ of degree (not greater than) $m$. The case $2l-\Nb=0,$ i.e., the critical one, was dealt with  in \cite{RSh2}, \cite{RIntegral}, and therefore we do not discuss it here.
  For a detailed explanation  of such  symbolic structure see, e.g.,  \cite{Taylor 2}, Ch. 2, especially, Proposition 2.6.
   It is possible that only the logarithmic term is present in the leading singularity of the kernel. This happens, e.g.,  when $\AF$ is, in the leading term, the fundamental solution of a power of the Laplacian, say, $\AF=(1-\Delta)^{-l/2}$ with positive even integer $2l-\Nb$. The structure of the integral kernel of \emph{this} fundamental solution has  been known since long ago, see, e.g., \cite{GelShil}, Ch.II.

After framing by $\G_{\Si}$ and $\G_{\Si}^*$, as in \eqref{transposed}, we arrive at the representation of $\KF\KF^*$ as the integral operator  $\RF$ in $ L_{2,\m}=L_{2,\m}(\Si)$ with kernel $R(X,Y,X-Y)$.  Exactly this kind of operators  was considered in the papers \cite{AgrAm} (where 'almost smooth' Lipschitz surfaces were studied)  and (in the general Lipschitz case) in \cite{RT1},  for  surfaces of codimension 1, and \cite{RT2}, for an arbitrary codimension. The result on the eigenvalue asymptotics, obtained for such integral operators in
\cite{RT1}, \cite{RT2}, corresponds exactly the formulas in Theorem \ref{Th.As.RSh} above. The correspondence  between the symbol of pseudodifferential operator $\KF\KF^*$ and its integral kernel is used in \cite{RT1}, \cite{RT2} systematically. Expression \eqref{r_-q} for the auxiliary symbol $\rb_{-m}$ is obtained by representing the restriction operator $\G_{\Si}$ by means of the Fourier transform.

  In \cite{RT1}, \cite{RT2}, the eigenvalue asymptotics formulas have been proved first for this operator on a \emph{smooth} surface $\Si$, passing again to a pseudodifferential representation of the integral operator. In this case, $\RF$ is  a classical pseudodifferential operator on $\Si,$ of order $-m,$ with principal symbol $\rb_{-m}$ given by  \eqref{r_-q}, expressed in the co-ordinate system generated by the orthogonal projection to the tangent plane at the point $X\in\Si$.  The proof of eigenvalue asymptotics now follows immediately from the, now classical, results by Birman-Solomyak on eigenvalue asymptotics for negative order pseudodifferential operators, see \cite{BirSolPDO1}. Note that the coefficient in the asymptotic formula has the meaning of phase volume.

  If the surface $\Si$ is not better than Lipschitz, $\RF$ is, generally, not a classical pseudodifferential operator, therefore, $\rb_{-m}(X,\x)$ is not a symbol of anything but it is just considered as an expression involved in calculating the density in the spectral asymptotics formula.
The given Lipschitz surface  $\Si:\yb=\pmb{\f}(\xb)$ is approximated, locally,
  by smooth ones, $\Si_\e,$ so that in their local representation $\yb=\pmb{\f}_\e(\xb)$,
  functions $\pmb{\f}_\e$ converge to $\pmb{\f}$ in $L_{\infty}$ and their gradients $\nabla\pmb{\f}_\e$ converge to $\nabla\pmb{\f}$ in all $L_p, p<\infty$ (one should not expect convergence of gradients in $L_\infty,$ of course). Expressed in local variables $\xb\mapsto (\xb,\pmb{\f}(\xb))$, resp., $\xb\mapsto (\xb,\pmb{\f}_{\e}(\xb))$,  operators with kernel $R(X,Y,X-Y)$ on surfaces $\Si$ and $ \Si_\e$ are transformed to operators $\tilde{\RF},$ resp., $\tilde{\RF}_{\e},$ on one and the some domain in $\R^{d}$, while the eigenvalue asymptotics for  ${\tilde{\RF}}_{\e},$ is known. Now it is possible to consider the difference  ${\tilde{\RF}}-{\tilde{\RF}}_{\e}$ of these operators. The eigenvalues of this difference are estimated using the closeness or $\pmb{\f}$ and $\pmb{\f}_\e$ as $\Si_\e\to\Si$ (and this is the most technical part of the reasoning in \cite{RT1}, \cite{RT2}),  so we obtain that the eigenvalue asymptotic coefficients of $ {\tilde{\RF}} -{\tilde{\RF}}_{\e}$ converge to zero. This property enables one to use again the asymptotic perturbation  Lemma \ref{BSLemma}, to justify the eigenvalue asymptotics formula for ${\RF}.$
\end{proof}
\subsection{Eigenvalue asymptotics on uniformly rectifiable sets.}\label{rect.sect}
(Uniformly) rectifiable sets are an important object of study in the geometric measure theory.
\begin{defin}\label{Def.URS} The set $\Ec\subset\R^\Nb$ is called uniformly rectifiable
 of dimension $d$ if $\Ec$ is the union
 of a countable collection of Lipschitz surfaces $\Si_j$ of dimension $d,$ up to a set of Hausdorff measure $0$, $\Hc^d(\Ec\setminus\cup \Si_j)=0.$
\end{defin}
A number of criteria for a set to be rectifiable has been found, see, e.g., a review in \cite{De Lellis}; a brief exposition with references is included
in  \cite{RIntegral}.  Note that any compact connected
set of Hausdorff dimension $1$ is uniformly rectifiable. In \cite{RIntegral}, in the critical case $2l=\Nb$, the result on eigenvalue asymptotics was proved for
$\m$ being  the union of uniformly rectifiable sets of different Hausdorff dimensions.

In the noncritical case the eigenvalue asymptotics is established in a similar way. We restrict ourselves to the case $\AF=\AF_0.$
\begin{thm}\label{symm.as}Let $\Mc$ be a compact uniformly rectifiable set of Hausdorff dimension $d$, $\max(0,2l-\Nb)<d<\Nb,$ $\AF_0$ be a pseudodifferential operator with principal symbol $|\X|^{-l}.$
Suppose that the Hausdorff measure $\mu=\Hc^d$ on $\Mc$ and the density $V$ on $\Mc$  satisfy Condition \ref{condition}. Then for  operator $\Tb_{V\m,\AF_0}$ the eigenvalue asymptotic formula
 \eqref{integral} holds, with integration over $\Mc.$
\end{thm}
\begin{proof}{The proof} repeats  the reasoning in  \cite{RIntegral} with minor changes. We explain here the main steps. First of all,
 as before, it is sufficient to consider the case of a sign-definite $V,$ e.g., $V\ge0$.
Let $\Si_j$, $j=1,\dots$ be a sequence of Lipschitz surfaces exhausting $\Mc$ up to a set of zero
Hausdorff measure.  Namely,  denote by $\Ec_k$ the finite union $\Ec_k=\cup_{j\le k}\Si_j$, $\Fc_k=\Mc\setminus\Ec_k$, so that $\Hc^d(\Fc_k)\to 0$ as
$k\to\infty.$ We split operator $\Tb_{V\m}=\Tb_{V\m,\AF_0}$ into the sum
\begin{equation}\label{small measure}
    \Tb_{V\m}=\Tb_{V_k\m}+\Tb_{V_k'\m}\equiv \Tb_k+\Tb_k',
\end{equation}
where $V_k$ is the restriction of  $V$ to the set $\Ec_k$ and $V_k'$ is the restriction of $V$ to $\Fc_k$.
 Since $\Hc^d(\Fc_k)$ tends to zero as $k\to\infty,$
 we have  $\int (V_k')^{\vartheta}\m(dX)\to 0$ as $k\to \infty.$ By estimates in Section \ref{Sect.3},
 the second term in  \eqref{small measure} satisfies $\Db^\theta_+(\Tb_k')\to 0$ as $k\to\infty.$ Therefore $\Tb_k'$ can serve as the second term in the
 decomposition in Lemma \ref{BSLemma}.

  It remains to prove the required asymptotic formula for the first
 term in \eqref{small measure}, i.e., for an operator with measure supported on a finite union of Lipschitz surfaces.
 This is performed by means of induction on $k,$ this means on the number of
 Lipschitz surfaces $\Si_j,$ $j\le k$. For $k=1,$ i.e., for one surface, the asymptotic formulas has been established
 in Theorem \ref{Th.As.RSh}.  Suppose now that the eigenvalue asymptotics has been already proved for the union of $k-1$ Lipschitz surfaces, i.e., for the compact set $\Ec_{k-1}.$ Let us add
  one more Lipschitz surface $\Si_{k}.$  Consider the set $\Nc_\de, $ the $\de$-neighborhood of $\Ec_{k-1}$
  and denote by $\Qc_{\de}$ the set $(\Nc_\de\cap \Si_k)\setminus \Ec_{k-1}.$ Function $V$ splits into three parts $V=V_0+V_\de+V_k,$ where $V_0$ is supported
  in $\Ec_{k-1},$  $V_k$ is supported in $\Si_{k}\setminus\Qc_\de,$ i.e., in just one Lipschitz surface, and, finally, $V_{\de}$ is supported in $\Qc_{\de}.$  Correspondingly, operator $T_{V,\m,\Mc}$ splits into the sum
  \begin{equation}\label{split.rectif}
    \Tb_{V\m} = (\Tb_0+\Tb_k)+\Tb_\de\equiv (\Tb_{V_0\m}+\Tb_{V_k\m})+\Tb_{V_\de\mu}.
  \end{equation}
  As $\de\to 0,$ the Hausdorff measure of $\Qc_{\de}$ tends to zero, and thus the $L_{\vartheta}$ -norm of $V_{\de}$ tends to zero.
  Therefore, by estimates in Sect.\ref{Sect.3}, $\Db_+^{\theta_0}(\Tb_\delta)$ tends
  to zero as $\de\to 0.$ For operators $\Tb_{V_0\m}$, $\Tb_{V_k\m}$ the asymptotic formulas hold, by the inductive assumption. Moreover,
 the supports of functions $V_0,V_k$ are well separated, their distance is not less than $\de$. This allows us to use Corollary \ref{Cor.addit}.
 to obtain the eigenvalue asymptotics for  $\Tb_0+\Tb_k$. Finally, the decomposition \eqref{split.rectif} enables us to apply, again,
 Lemma \ref{BSLemma}.
 \end{proof}
 \subsection{Eigenvalue asymptotics on noncompact locally Lipschitz sets.}\label{Sect.noncompact} It is interesting to extend the results on eigenvalue estimates and asymptotics to a reasonably large class of noncompact Lipschitz surfaces or sets composed of these. Generally, due to non-compactness, some additional conditions are needed to justify eigenvalue properties in question.

  We call a set $\Si\subset\R^{\Nb}$ with Hausdorff measure $\m=\Hc^d|_\Si$  a locally Lipschitz set if there exists a locally finite family  of  compact Lipschitz surfaces $\Si_j$ of dimension $d$ such that $\Si=\cup \Si_j$
 \begin{rem} Any finite collection of surfaces $\Si_j, j\le k,$  has a common Lipschitz constant. For the infinite set of surfaces, $\Si_j, j<\infty,$ such common constant does not necessarily exist, so the surfaces may become more and more curved as one goes to infinity.
 \end{rem}
 \begin{thm}\label{th.lip}Let $2l<\Nb$ and $\Si$ be a locally Lipschitz  set of dimension $d$ in $\R^\Nb,$ $\AF_0=(-\Delta)^{-l/2},$  $2l<\Nb, \,d>\Nb-2l.$ Let the Hausdorff measure $\m$ on $\Si$ satisfies \eqref{AR-} with $s=d$. Suppose that $V\in L_{\theta,\mu}(\Si).$ Then for  operator $\Tb=\Tb_{P,\AF_0}$, $P=V\mu$
 the eigenvalue asymptotic formula \eqref{AsLambda}-\eqref{IntCoeff} holds.
 \end{thm}
  \begin{proof} The proof follows the pattern of the one for Theorem \ref{symm.as}. For a given $\e>0,$ we split $V$ into two parts, $V=V_\e+V_\e'$ so that $V_\e$ is supported in the union of a finite set of surfaces $\Si_j,$ while $V_\e'$ has small $L_{\theta,\m}$-norm, $\|V\|_{L_{\theta,\m}}<\e.$ For operator $\Tb_{P_\e,\AF_0},$ $P_\e=V_\e\m$ the eigenvalue asymptotic formula was established when proving Theorem \ref{symm.as}. For operator $\Tb_{P_\e',\AF_0}$, $P_\e'=V_\e'\m,$
we have, by Theorem \ref{ThmPsdoRN} eigenvalue estimate with small constant. Now, as usual, Lemma \ref{BSLemma} finishes the job.   \end{proof}

For the supercritical case, we restrict ourselves to the Birman-Borzov type of operators.

\begin{thm}\label{Lip.BB} Let $\Si$ is a locally Lipschitz  set of dimension $d$ in $\R^\Nb,$ $\AF_0=(-\Delta)^{-l/2},$  $2l>\Nb$. Suppose that Hausdorff measure $\m$ satisfies \eqref{AR-noncomp}. Let $V$ be a real function on $\Si$ satisfying \eqref{BBcond}. Then for the eigenvalues of operator $\Tb=\Tb_{P,\AF_0}$, $P=V\mu$
 the eigenvalue asymptotic formula \eqref{AsLambda}-\eqref{IntCoeff} holds.
\end{thm}
\begin{proof}The reasoning is the same as in Theorem \ref{th.lip}. We just apply Corollary \ref{BB.cor} to estimate the eigenvalues
outside a compact set.\end{proof}

\appendix
\section{The approach by D.Edmunds and H.Triebel}
 In this Appendix we discuss relation of our results with the ones presented in books by D.Edmunds and H.Triebel, see \cite{ET}, and further by H.Triebel, \cite{Triebel1}, \cite{Triebel2}, H.Triebel and D.Haroske, \cite{HT}, and  accompanying papers.

 Starting from early papers by H.Triebel, an approach  for obtaining quantitative characteristics of operators in Banach and quasi-Banach spaces  was being developed, based upon the analysis of entropy numbers. We recall that for a \emph{compact } linear operator $\Hb: \Xc\to\Yc,$  the entropy number $e_k(\Hb)$ is defined as the smallest $\e$ such that the image in $\Yc$ of the unit ball in $\Xc$ can be covered by not more than $2^k$ balls of radius $\e$ in the metric of $\Yc.$  The approximation numbers $a_k(\Hb)$ of a compact operator $\Hb$ in a quasi-Banach space $\Xc$ are defined as
 \begin{equation*}
    a_k(\Hb)=\inf\{\|\Hb-\Vb\|_{\Xc\to\Yc}:\rank(\Vb)\le k.\}
 \end{equation*}
When $\Xc=\Yc,$ the approximation numbers are closely related to the eigenvalues $\m_k(\Hb),$ and for the case of Hilbert spaces, $\m_k(\Hb)$ are the singular numbers of $\Hb.$

The approach we are discussing now consists in the following. Having an operator $\Hb$ in a Hilbert or Banach space, of a complicated structure, containing multiplications by weight functions, trace  and cotrace operators, one factorizes $\Hb$ as a composition of several operators, using embedding and trace theorems for functional and distributional spaces and H\"older type inequalities, so that all but one operators in this composition, the ones containing weight functions, are bounded, and there is one embedding  operator in some standard spaces, which is compact and for which estimates for entropy numbers are known. This implies estimates for entropy numbers for $\Hb.$ Finally, estimates for approximation numbers or eigenvalues (singular numbers) follow from the wonderful B.Carl's theorem relating these quantities.

The most simple example of this construction can be seen in \cite{ET}, Sect. 5.2.4. The authors consider there a Birman-Schwinger type operator
\begin{equation}\label{App.1}
    \Hb = b_2 C b_1
\end{equation}
in $L_2(\Om),$ where $C=A^{-1},$ and $A$ is an order $2m$ elliptic operator in a bounded domain $\Om\subset\R^{\Nb}$ with some elliptic boundary conditions so that $C$ is defined and maps Sobolev spaces $H_q^l$ with gain of $2m$ derivatives, $C: H_q^l(\Om)\to H_q^{l+2m}(\Om),$ $q\in[1,\infty).$ The weight functions $b_1, b_2$ belong, respectively, to the spaces $L_{r_1}(\Om), L_{r_2}(\Om).$ By Proposition 5.2.4, (we present it in the Hilbert space setting, $p=2$) if $r_1>2,r_2>2,$ $\de=\frac{2m}{\Nb}-r_1^{-1}-r_2^{-1}>0,$ operator $b_2 A^{-1}b_1$ is compact and its eigenvalues in $L_2(\Om)$ satisfy
\begin{equation}\label{App.2}
    |\la_k(\Hb)|\le C(A,\Om) ||b_1||_{L_{r_1}} ||b_2||_{L_{r_2}}k^{-\frac{2m}{\Nb}}.
\end{equation}

The factorization mentioned above has the form
\begin{equation}\label{3}
    \Hb=b_2\circ \mathrm{id}_{H_q^{2m}\to L_t} \circ A^{-1}\circ b_1,
\end{equation}
where the factors are: $b_1: L_2\to L_q,$ $q^{-1}=r_1^{-1}+1/2,$ $A^{-1}: L_q\to H^{2m}_q,$ $\mathrm{id}: H_q^{2m}\to L_t,\, t^{-1}=1/2-r_2^{-1}, $ $b_2:L_t\to L_2.$ Here, the operator $\mathrm{id},$ the embedding of the Sobolev space $H_q^{2m}$ into  $L_t$ is compact and an estimate of its entropy numbers was previously found, $e_k(id)\le C k^{-\frac{2m}{\Nb}}$.
This implies a similar estimate for entropy numbers of $\Hb$, and, by the Carl theorem, produces the required estimate for singular values of $\Hb$.

The  circumstance of great importance here is the fact that the entropy numbers $e_k(\mathrm{id})$ of the embedding  $\mathrm{id}: H_q^{2m}\to L_t$ have decay rate not depending on the parameter $t,$ in other words, on the space into which the Sobolev space is embedded, as long as this embedding is compact, $t^{-1}>q^{-1}-\frac{2m}{\Nb}$, equivalently, $\de>0$.  This looks mysterious if one compares this property of entropy numbers with the corresponding properties of approximation numbers of the same embedding: the latter decay rate \emph{does} depend on the parameter $t.$ This property of entropy numbers was first discovered by M.Birman and M.Solomyak in \cite{BSPiecewise}, and it is these results that formed the basis  of considerations in \cite{ET}. Further, these estimates for entropy numbers have been carried over, mostly by H.Triebel, to quite a lot of embedding operators in  various spaces, using interpolation and other  advanced constructions. About this property of entropy numbers, the authors of \cite{HT} write: "This somewhat
surprising assertion is a consequence of the miraculous properties of entropy numbers...", see P. 243.

What is inherent and unavoidable  for this approach is that the summability conditions for the weight functions, such as $b_1,b_2$ in the above example, are not sharp.  Say, in the most classical case $b_1=b_2, $ $r_1=r_2=r$, $\Nb>2m$ (the subcritical case in our meaning), these conditions require that $r>\Nb/m.$ If $\de$ in the above example equals zero, the embedding operator $\mathrm{id}_{H_q^{2m}\to L_t}$ is continuous but not compact and the entropy numbers estimates fail. In the opposite, estimates of CLR type admit $r=m/\Nb,$ equivalently, $\de=0$, see, e.g. estimates presented in   \cite{BS} and further publications, up to this one.  This, seemingly minor, difference is, in fact, quite crucial: the constants in such estimates with sharp exponents do not depend on the size of the domain $\Om$ and therefore eigenvalue estimates can be carried over to operators in unbounded domains, in particular, to the whole space, say, as it was done in Theorem \ref{ThmPsdoRN}. It is this non-sharpness of results of \cite{ET} that leads to excessive conditions in estimates of the negative spectrum of Schr\"odinger like operators, namely, the compactness of the support of the potential, see, e.g.,  \cite{ET}, Sect.5.2.7.

A similar effect is observed for $\Nb<2m$ (the supercritical case). The condition $r_j>2, $ caused, again, by the requirement that a certain embedding operator is compact, is not sharp. On the other hand,  in the M. Birman--M. Solomyak approach, in a similar setting, it is allowed that $r_j=2.$

These limitations of their method for obtaining eigenvalue estimates was well understood by the authors of \cite{ET}: they write, see P.185, 186:
"As would be expected, in these special situations, the deep Hilbert space
techniques used by the authors listed above  [M.Birman, M.Solomyak, B.Simon, G.Rozenblum, G.Tashchiyan] often give better results
than those obtained by our simple arguments which are not confined to
symmetric operators or Hilbert spaces."

Similar reservations are made in \cite{Triebel2}. Say on  P.231, where complications arising in the 'limiting case' are discussed, in particular, as it concerns the passage to singular measures with unbounded support,
it is emphatically  stated: 'But we do not go into detail for a simple
reason: Nothing has been done so far.'

The above approach has been applied to various types of spectral problems. In particular, in \cite{Triebel2}, operators of the form $B=b_2(\g)b(x,D)b_1(\g)$ have been considered, with $b_j$ being functions on a set $\g\subset \R^{\Nb},$ under assumption that they belong to $L_{r_j}(\g)$ with respect to the Hausdorff measures of dimension $d,$ while the set $\g$ is a $d$-set (is Ahlfors regular of order $d,$ $0<d<\Nb$). Further on, $b(x,D)$ is an order $-\varkappa$ pseudodifferential operator and a number of conditions on parameters of the problem are imposed, see Sect. 27, 28 in \cite{Triebel2}. For obtaining eigenvalue estimates, again, operator under study is factorized as a composition of several operators, including multiplication by weight functions, restriction to $\g$, extension from $\g$, and embedding operators for functional spaces on $\R^\Nb$ and on $\g.$ Deep and complicated results of this book and previous ones are used to justify boundedness of these operators. And for one of them, which turns out to be compact, namely, the embedding of Besov type spaces, the estimate of entropy numbers has been established.

As before,  the compactness requirement for this particular embedding operator leads to non-sharp integrability conditions for weight functions $b_1,b_2,$ with the same kind of shortcomings of the resulting eigenvalue estimates for operator $B.$

In Sect.\ref{Sect.4} in our paper, we consider some of operators of the type studied in \cite{Triebel2} and demonstrate that our approach produces spectral estimates, sharp both in order and integrability classes of weight functions.

So, when comparing the results by H.Triebel and the ones of the present paper, we can say the following.
\begin{enumerate}
\item Results by H.Triebel and co-authors treat a wide class of operators, including the ones close to operators in the present paper, and not necessarily in the Hilbert space setting. The integrability condition imposed upon the weight functions are not sharp, so the estimates lack the homogeneity property for weight functions, needed for getting rid of the compact support conditions. Additionally, the two-sided Ahlfors condition \eqref{AR} is required.
\item Results in our paper deal with a more narrow class of operators, exclusively in the Hilbert space setting. However, the integrability conditions on the weight functions are sharp,and their homogeneity properties enable us to get rid of the support compactness conditions, in particular, producing CLR-type eigenvalue estimates. Only one-sided Ahlfors conditions \eqref{AR+} or \eqref{AR-} are imposed on the measure $\m$.
\end{enumerate}

\end{document}